\documentclass[11pt]{article}
\usepackage[english]{babel}
\usepackage[latin1]{inputenc}
\usepackage{amsmath}
\usepackage{amssymb}
\usepackage{amsfonts}
\usepackage{amsthm}
\usepackage{bbm}
\usepackage{color}
\usepackage{graphicx}

\newcommand{\nm}{\nabla_{X,z}^{\mu_2}}
\newcommand{\dm}{\Delta_{X,z}^{\mu_2}}
\newcommand{\cst}{\mbox{\textnormal{Cst }}}
\newcommand{\dsp}{\displaystyle}
\newcommand{\dt}{\partial_t}
\newcommand{\dz}{\partial_z}
\newcommand{\dx}{\partial_x}

\newcommand{\R}{{\mathbb R}}
\newcommand{\N}{{\mathbb N}}
\newcommand{\eps}{\varepsilon}

\newcommand{\cG}{{\mathcal G}[\zeta]}
\newcommand{\cV}{{\mathcal V}[\zeta]}
\newcommand{\cZ}{{\mathcal Z}[\zeta]}
\newcommand{\tm}{{\bf T}_{\mu}}
\newcommand{\tmd}{{\bf T}_{\mu_2}}
\newcommand{\D}{\vert D\vert}
\theoremstyle{remark}
\newtheorem{example}{Example}
\newtheorem{remark}{Remark}
\newtheorem{notation}{Notation}
\theoremstyle{definition}
\newtheorem{definition}{Definition}
\theoremstyle{theorem}
\newtheorem{proposition}{Proposition}
\newtheorem{lemma}{Lemma}
\newtheorem{theorem}{Theorem}
\newtheorem{corollary}{Corollary}

\begin{document}
\title{Asymptotic Models for Internal Waves}

 \author{
\renewcommand{\thefootnote}{\arabic{footnote}}
J. L. Bona\footnotemark[1], D. Lannes\footnotemark[2]~~and J.-C.
Saut\footnotemark[3]} \footnotetext[1]{Department of Mathematics,
Statistics and Computer Science, University of Illinois at Chicago,
Chicago, IL 60607, USA. E-mail: bona@math.uic.edu}
\footnotetext[2]{Universit\'e Bordeaux I; IMB et CNRS UMR 5251, 351
cours de la Lib\'eration, 33405 Talence Cedex, France. E-mail:
David.Lannes@math.u-bordeaux1.fr} \footnotetext[3]{Universit\'e de
Paris-Sud et CNRS UMR 8628, B\^at. 425, 91405 Orsay Cedex, France.
E-mail: jean-claude.saut@math.u-psud.fr}

\date{December 17, 2007}
\maketitle

\begin {abstract}
 We derived here in a systematic way, and for a large class of scaling
regimes, asymptotic models for the propagation of internal waves
at the interface between two layers of immiscible fluids of
different densities, under the rigid lid assumption and with a flat
bottom. The full (Euler) model for this situation is reduced to a
system of evolution equations posed spatially on $\R^d$, $d=1,2$,
which involve two nonlocal operators. The different asymptotic
models are obtained by expanding the nonlocal operators with respect
to suitable small parameters that depend variously on the amplitude,
wave-lengths and depth ratio of the two layers. We rigorously derive
classical models and also some model systems that appear to be new.
Furthermore, the consistency of these asymptotic systems with the
full Euler equations is established.\\

\vspace{1 cm}
Nous \' etablissons ici de mani\` ere syst\' ematique, et pour une grande classe de r\' egimes,
des mod\`eles asymptotiques pour la propagation d'ondes internes \`a l'interface de deux couches de fluides immiscibles de densit\' e diff\' erente, sous l'hypoth\`ese de toit rigide et de fond plat. Les \' equations compl\`etes pour cette situation (Euler) sont r\' eduites \`a un syst\`eme d'\' equations d'\' evolution pos\' e dans le domaine spatial $\R^d$, $d=1,2$, et qui comprend deux op\' erateurs non locaux. Les divers mod\`eles asymptotiques sont obtenus en d\' eveloppant les op\' erateurs non locaux par rapport \`a des petits param\`etres convenables (d\' ependant de l'amplitude, de la longueur d'onde et du rapport de hauteur des deux couches). Nous \' etablissons rigoureusement des mod\`eles classiques ainsi que d'autres qui semblent nouveaux. De plus, on montre la consistance de ces syst\`emes asymptotiques avec les \' equations d'Euler. 

\end {abstract}

\section{Introduction}

\subsection{General Setting}

The mathematical theory of waves on the interface between two layers
of immiscible fluid of different densities has attracted interest
because it is the simplest idealization for internal wave
propagation and because of the challenging modeling, mathematical
and numerical issues that arise in the analysis of this system. The
recent survey article of Helfrich and Melville \cite{HM} provides a
rather extensive bibliography and a good overview of the properties
of steady internal solitary waves in such systems as well as for
more general density stratifications. The compendium \cite{J} of
field observations comprised of synthetic aperture radar (SAR)
images of large-amplitude internal waves in different oceans
together with associated physical data shows just how varied can be
the propagation of internal waves.  This variety is reflected in the
mathematical models for such phenomena.  Because of the range of
scaling regimes that come to the fore in real environments, the
literature on internal wave models is markedly richer in terms of
different types of model equations than is the case for surface wave
propagation (see, e.g. \cite{BCS, BCL} and the references therein).

The idealized system that will be the focus of the discussion here,
when it is at rest, consists of a homogeneous fluid of depth $d_1$
and density $\rho_1$ lying over another homogeneous fluid of depth
$d_2$ and density $\rho_2 > \rho_1$.  The bottom on which both
fluids rest is presumed to be horizontal and featureless while the
top of fluid 1 is restricted by the rigid lid assumption, which is
to say, the top is viewed as an impenetrable, bounding surface. This
is a standard assumption, and is reckoned to be a good one when the
pycnocline is far from the top, which is when $d_1$ is large
relative to the wavelength of a disturbance. In the present work,
two general classes of waves will be countenanced.  Both of these
require that the deviation of the interface be a graph over the flat
bottom, so overturning waves are not within the purview of our
theory (see Figure 1 for a definition sketch). The first, which is
referred to as the one-dimensional case, are long-crested waves that
propagate principally along one axis, say along the $x$-direction in
a standard $x-y-z$ Cartesian frame in which $z$ is directed opposite
to the direction in which gravity acts.  Such motions are taken to
be sensibly independent of the $y$-coordinate and can be
successfully modeled in the first instance by the two-dimensional
Euler system involving only the independent variables $x$, $z$ and
of course time $t$. Because the interface is a graph over the
bottom, these asymptotic models then depend only upon $x \in \R$ and
$t$, and hence the appellation `one-dimensional'. Among
one-dimensional models, the simplest are those in which one further
assumes that the waves travel only in one direction, say in the
direction of increasing values of $x$. Models which we will call
`two-dimensional' are not restricted by the long-crested
presumption, and are consequently more general than the
one-dimensional models.  They are derived from the full
three-dimensional Euler system and their dependent variables depend
upon the spatial variable $X = (x,y) \in \R^2$ and time $t$.

One-dimensional, unidirectional, weakly nonlinear models such as the
Korteweg-de Vries (KdV) equation , the Intermediate Long Wave (ILW)
equation \cite{Joseph, KK} or the Benjamin-Ono equation \cite{B}
have been extensively used and compared with laboratory experiments
\cite{ KB, SH, W}.  While much of our qualitative appreciation of
the interaction between the competing effects of nonlinearity and
dispersion in surface and internal wave propagation has been
informed by these sorts of equations, they are of somewhat limited
validity (c.f. \cite{ABR}). Weakly nonlinear models in
two-dimensions have been derived by Camassa and Choi \cite{CC1}.
Nguyen and Dias \cite{ND} have derived and studied a Boussinesq-type
system in a weakly nonlinear regime. Fully nonlinear models were
obtained in the one-dimensional case by Matsuno \cite{Matsuno}, and
in the two-dimensional case by Camassa and Choi \cite{CC2}. We
mention also the interesting paper by Camassa {\it et al.}
\cite{CCMRS} where the aforementioned models are compared, in the
one-dimensional case, with experimental observations and numerical
integrations of the full Euler system. In \cite{CC1,CC2,Matsuno} the
analysis commences with the full Euler system formulation and the
asymptotic models are obtained by formally expanding the unknowns
with respect to a small parameter. It is not easy using this
approach to provide a rigorous justification of the asymptotic
expansion, except perhaps within the setting of analytic functions.
A different approach has been carried out by Craig, Guyenne and
Kalisch \cite{CGK} in the one-dimensional case. These authors use
the Hamiltonian formulation of the Euler equations (due originally
to Zakharov \cite{Zakharov} for surface waves and to Benjamin and Bridges 
\cite{BB} for internal waves) and expand the Hamiltonian with respect to the relevant
small parameters.  This method provides a hierarchy of Hamiltonian
systems which serve as approximations of the full Euler equations.
Such systems are not always the best for modeling, analytical or
numerical purposes, however.  Indeed, they can even be linearly
ill-posed in Hadamard's classical sense. In such cases, it is
necessary in the Hamiltonian framework to proceed one stage further
in the expansion, leading to more complicated systems (which may
still not be well posed).

The strategy followed here is inspired by that initiated in
\cite{BCS, BCS2, BCL}. Namely, following the procedure introduced in
\cite{CGK,CSS,Zakharov}, we rewrite the full system as a system of
two evolution equations posed on $\R^d$, where $d=1$ or 2 depending
upon whether a one- or two-dimensional model is being contemplated.
The reformulated system, which depends only upon the spatial
variable on the interface, involves two non-local operators, a
Dirichlet-to-Neumann operator $G[\zeta]$, and what we term an
``interface operator" ${\bf H}[\zeta]$, defined precisely below. Of
course the operator ${\bf H}[\zeta]$ does not appear in the theory
of surface waves, and this is an interesting new aspect of the
internal wave theory.  A rigorously justified asymptotic expansion
of the non-local operators with respect to dimensionless small
parameters is then mounted.  We consider both the ``weakly
nonlinear" case and the ``fully nonlinear" situation and cover a
variety of scaling regimes. For the considered scaling regimes,
these expansions then lead to an asymptotic evolution system.  As in
\cite{BCS,BCS2, BCL}, in each case a family of asymptotic models may
then be inferred by using the ``BBM trick" and suitable changes of
the dependent variables.  This analysis recovers most of the systems
which have been introduced in the literature and also some
interesting new ones.   For instance, in certain of the
two-dimensional regimes, a non-local operator appears whose analog
is not present in any of the one-dimensional cases.

All the systems derived are proved to be consistent with the full
Euler system. In rough terms, this means that any solution of the
latter solves any of the asymptotic systems up to a small error. The
systems are thus seen to be formally equivalent models in terms of
the small parameter`s that arise in the expansions.  The advantage
of obtaining a family of equivalent asymptotic systems is clear from
the modeling perspective.  One can use the flexibility to adjust the
linearized dispersion relation to better fit the exact dispersion
and can choose horizontal velocity variables that are well suited to
the predictions in view. Mathematically, the choice will be among
those that are well posed for the particular initial-value or
initial-boundary-value problem under consideration. When it comes to
computer simulation, some of the systems are far better suited to
the construction of stable, accurate numerical schemes and these
would naturally be favored.

The paper is organized as follows. In the next portion of the
Introduction, the ``Zakharov formulation" of the full system is
written in dimensionless form and the different scaling regimes
which will be studied enunciated.  In Subsection 1.5, a compendium
of the outcome of our analysis is offered to guide the reader
through the rest of the paper.   Chapter 2 is devoted to the
rigorous asymptotic analysis of the non-local, Dirichlet-to-Neumann
operator $G[\zeta]$ and the interface operator ${\bf H}[\zeta]$
mentioned earlier.  The asymptotic models that result from the use
of the expansions of these two operators are introduced (and proved
to be consistent with the full Euler system) in Chapter 3. The
somewhat technical proof of Proposition 3 is given in Appendix A.

In the present paper, we have refrained from pursuing the analysis
to the point of obtaining convergence results for the asymptotic
systems to the full internal waves system.  Such a program has been
fully achieved in the case of surface waves by combining the results
of \cite{AL} and \cite{BCL}.  What is needed to complete the circle
of ideas in the internal wave case is a stability analysis  of the
asymptotic models derived here (that is, an estimation of the
remainders which comprise the difference between the Euler system
and the models). Together with consistency, a straightforward
analysis would then provide a convergence result to the full Euler
system, assuming that the large time existence results obtained by
Alvarez-Samaniego and Lannes in \cite{AL} for the surface wave
system are valid for the internal waves system. The latter point is
far from obvious; indeed, even the local well-posedness of the Euler
equations in the two-fluid configuration seems to be an open problem
in the absence of surface tension (cf \cite{OI} for the rigorous derivation
of the Benjamin-Ono equation for the two-fluid system in the presence of
 surface tension).

Note finally that the analysis of the present paper could be
extended to the case of a seabed with structure (a non-flat bottom,
see \cite{Chazel} for the case of surface waves) and to the case of
a two-layer system where the upper surface is free rather than
restricted by the rigid lid hypothesis (see \cite{BGT,CSS,Matsuno}
for a derivation of asymptotic models in this situation).  These
issues are under study and an analysis will be reported separately.
Of especial interest is a comparison of the problem considered with
the rigid-lid condition at the top and the problem wherein the upper
surface is left free in the case where $d_1$ is relatively large.

\begin{figure}[!b]
    \begin{centering}
    \includegraphics[width=12cm]{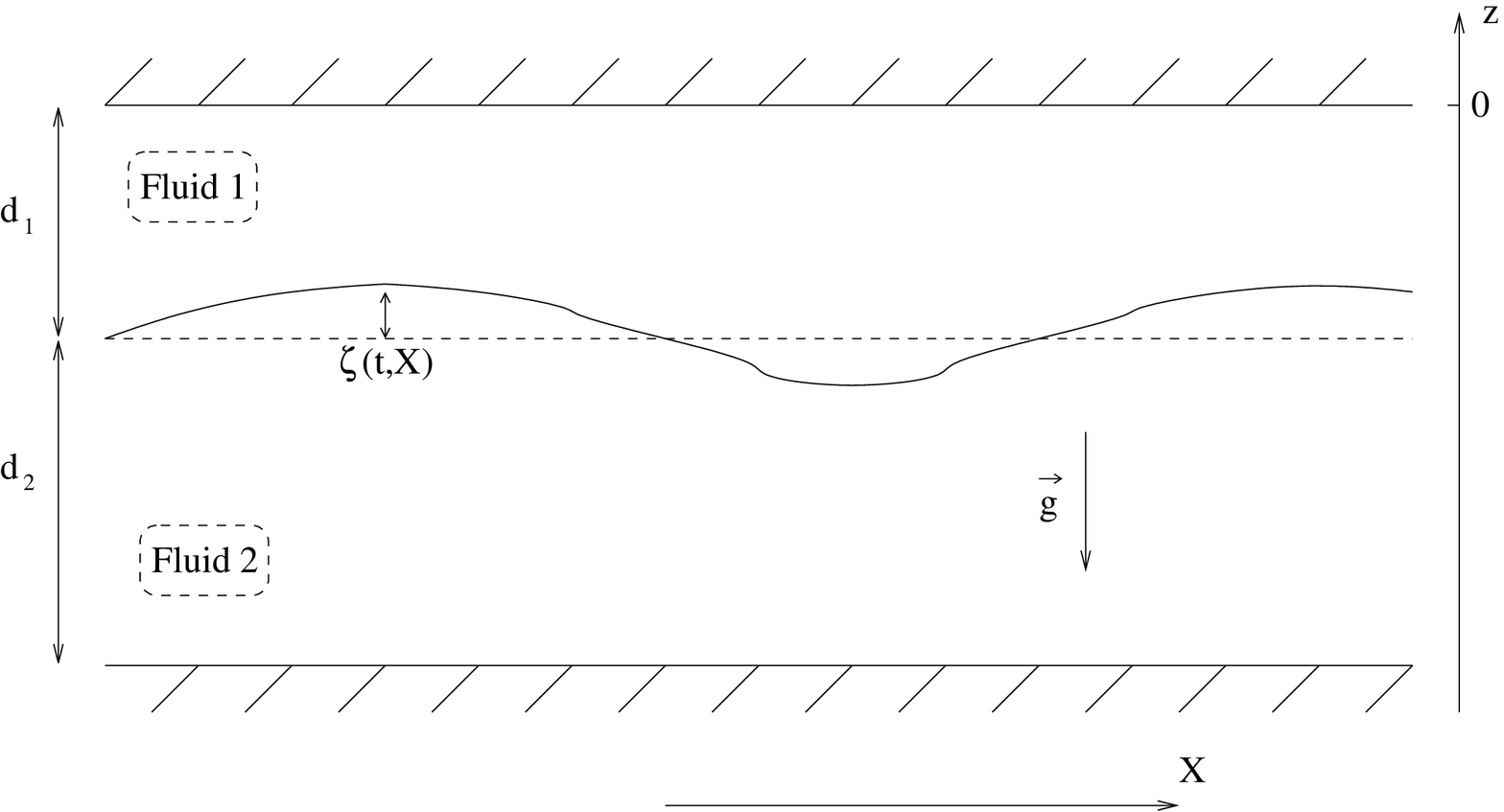}
    \end{centering}
\end{figure}
\clearpage

\centerline{\it Notation}

 Denote by $X$ the $d$-dimensional
horizontal variable as described earlier, where $d=1,2$.  Thus, $X =
x$ when $d=1$ and $X = (x,y)$ when $d=2$.  We continue to use $z$
for the vertical variable.

The usual symbols $\nabla$ and $\Delta$ connote the gradient and
Laplace operator in the horizontal variables, whereas $\nabla_{X,z}$
and $\Delta_{X,z}$ are their $d+1$-variable version (the gradient in
both or all three variables, depending on whether $d= 1$ or 2 and
similarly for the Laplacian).  For $\mu > 0$, it is very convenient
to also introduce scaled versions of the gradient and Laplace
operators, namely $\nabla_{X,z}^\mu=(\sqrt{\mu}\nabla^T,\dz)^T$ and
$\Delta_{X,z}^\mu=\nabla_{X,z}^\mu\cdot\nabla_{X,z}^\mu=\mu\Delta+\dz^2$.

For any tempered distribution $u$, denote by $\widehat{u}$ or
${\mathcal F}u$ its Fourier transform. If $f$ and $u$ are two
functions defined on $\R^d$, we use the Fourier multiplier notation
$f(D)u$ which is defined in terms of Fourier transforms, {\it viz.}
$$
    \widehat{f(D)u}=f\widehat{u}.
$$

The projection onto gradient fields in $L^2(\R^d)^d$ is written
$\Pi$ and is defined by the formula
$$
    \Pi=-\frac{\nabla \nabla^T}{\vert D\vert^2}.
$$
(Note that $\Pi=Id$ when $d=1$.)  The operator
$\Lambda=(1-\Delta)^{1/2}$ is equivalently defined using the Fourier
multiplier notation to be $\Lambda=(1+\vert D\vert^2)^{1/2}$.
Appearing frequently are the Fourier multipliers $\tm$ and $\tmd$,
given by
$$
    \tm=\tanh(\sqrt{\mu}\D)\quad\mbox{ and }\quad
    \tmd=\tanh(\sqrt{\mu}_2\D);
$$
where $\mu, \mu_2 > 0$.

The standard notation $H^s(\R^d)$, or simply $H^s$ if the underlying
domain is clear from the context, is used for the $L^2$-based
Sobolev spaces; their norm is written $\vert \cdot\vert_{H^s}$.

The planar strip ${\mathcal S}=\R^d\times (-1,0)$ appears often. The
unadorned norm $\Vert \cdot\Vert$ will always be the usual norm of
$L^2({\mathcal S})$.

\subsection{The Equations}

The Euler system of equations for our system is reviewed here.  As
in Figure 1, the origin of the vertical coordinate $z$ is taken at
the rigid top of the two-fluid system. Assuming each fluid is
incompressible and each flow irrotational, there exists velocity
potentials $\Phi_i$ ($i=1,2$) associated to both the upper and lower
fluid layers which satisfy
\begin{equation}
    \label{laplace}
    \Delta_{X,z}\Phi_i=0\quad \mbox{ in } \Omega^i_t
\end{equation}
for all time $t$, where $\Omega^i_t$ denotes the region occupied by
fluid $i$ at time $t$, $i=1,2$.  As above, fluid 1 refers to the
upper fluid layer whilst fluid 2 is the lower layer (see again
Figure 1). Assuming that the densities $\rho_i$, $i=1,2$, of both
fluids are constant, we also have two Bernouilli equations, namely,
\begin{equation}
    \label{bernouilli}
    \dt \Phi_i+\frac{1}{2}\vert \nabla_{X,z}\Phi_i\vert^2=-\frac{P}{\rho_i}-gz
    \quad \mbox{ in } \Omega^i_t,
\end{equation}
where $g$ denotes the acceleration of gravity and $P$ the pressure
inside the fluid. These equations are complemented by two boundary
conditions stating that the velocity must be horizontal at the two
rigid surfaces $\Gamma_1:=\{z=0\}$ and $\Gamma_2:=\{z=-d_1-d_2\}$,
which is to say
\begin{equation}
    \label{mur}
    \dz \Phi_i=0 \quad \mbox{ on }\quad \Gamma_i,\qquad
    (i=1,2).
\end{equation}
Finally, as mentioned earlier, it is presumed that the interface is
given as the graph of a function $\zeta(t,X)$ which expresses the
deviation of the interface from its rest position $(X,-d_1)$ at the
spatial coordinate $X$ at time $t$. The interface
$\Gamma_t:=\{z=-d_1+\zeta(t,X)\}$ between the fluids is taken to be
a bounding surface, or equivalently it is assumed that no fluid
particle crosses the interface. This condition, written for fluid
$i$, is classically expressed by the relation $\dt
\zeta=\sqrt{1+\vert \nabla\zeta\vert^2}v_n^i$, where $v_n^i$ denotes
the upwards normal derivative of the velocity of fluid $i$ at the
surface. Since this equation must of course be independant of which
fluid is being contemplated, it follows that the normal component of
the velocity is continuous at the interface. The two equations
\begin{equation}
    \label{kinematic}
    \dt \zeta=\sqrt{1+\vert\nabla\zeta\vert^2}\partial_n\Phi_1\quad \mbox{ on }\quad \Gamma_t,
\end{equation}
and
\begin{equation}
    \label{continuity}
    \partial_n\Phi_1=\partial_n\Phi_2 \quad \mbox{ on }\quad \Gamma_t,
\end{equation}
with
$$
    \partial_n:={\bf n}\cdot\nabla_{X,z}
    \quad \mbox{ and }\quad
    {\bf n}:=\frac{1}{\sqrt{1+\vert\nabla \zeta\vert^2}}
    (-\nabla \zeta,1)^T
$$
follow as a consequence.  A final condition is needed on the
pressure to close this set of equations, namely,
\begin{equation}
    \label{pression}
    P\mbox{ is continuous at the interface}.
\end{equation}

\subsection{Transformation of the Equations}

In this subsection, a new set of equations is deduced from the
internal-wave equations (\ref{laplace})-(\ref{pression}). Introduce
the trace of the potentials $\Phi_1$ and $\Phi_2$ at the interface,
$$
    \psi_i(t,X):=\Phi_i(t,X,-d_1+\zeta(t,X)), \qquad (i=1,2).
$$
One can evaluate Eq. (\ref{bernouilli}) at the interface and use
(\ref{kinematic}) and (\ref{continuity}) to obtain a set of
equations coupling $\zeta$ to $\psi_i$ ($i=1,2$), namely
\begin{eqnarray}
    \label{eq1}
    \dt \zeta-\sqrt{1+\vert\nabla\zeta\vert^2}\partial_n\Phi_i &=0, \\
        \label{eq2}
    \rho_i\Big(\dt \psi_i+g\zeta+\frac{1}{2}\vert \nabla \psi_i\vert^2
    -\frac{(\sqrt{1+\vert\nabla\zeta\vert^2}(\partial_n\Phi_i)+
    \nabla\zeta\cdot\nabla\psi_i)^2 }{2(1+\vert
    \nabla\zeta\vert^2)}\Big)
    &=-P,
\end{eqnarray}
where in (\ref{eq1}) and (\ref{eq2}), $(\partial_n\Phi_i)$ and $P$
are both evaluated at the interface $z=-d_1+\zeta(t,X)$.  Notice
that $\partial_n\Phi_1$  is fully
 determined by $\psi_1$
since $\Phi_1$ is uniquely given as the solution of Laplace's
equation (\ref{laplace}) in the upper fluid domain, the Neumann
condition (\ref{mur}) on $\Gamma_1$  and the Dirichlet condition
$\Phi_1=\psi_1$ at the interface. Following the formalism introduced
for the study of surface water waves in \cite{CSS2,CSS,Zakharov}, we
can therefore define the Dirichlet-Neumann operator $G[\zeta]\cdot$
by
$$
    G[\zeta]\psi_1=
    \sqrt{1+\vert\nabla\zeta\vert^2}
    (\partial_n\Phi_1)_{\vert_{z=-d_1+\zeta}}.
$$
Similarly, one remarks that $\psi_2$ is determined up to a constant
by $\psi_1$ since $\Phi_2$ is given (up to a constant) by the
resolution of the Laplace equation (\ref{laplace}) in the lower
fluid domain, with Neumann boundary conditions (\ref{mur}) on
$\Gamma_2$ and $\partial_n\Phi_2=\partial_n\Phi_1$ at the interface
(this latter being provided by (\ref{continuity})). It follows that
$\psi_1$ fully determines $\nabla\psi_2$ and we may thus define the
operator ${\bf H}[\zeta]\cdot$ by
$$
    {\bf H}[\zeta]\psi_1=\nabla\psi_2.
$$

Using the continuity of the pressure at the interface expressed in
(\ref{pression}), we may equate the left-hand sides of
(\ref{eq2})$_1$ and (\ref{eq2})$_2$ using the operators $G[\zeta]$
and ${\bf H}[\zeta]$ just defined. This yields the equation
$$
    \dt(\psi_2-\gamma\psi_1)+g(1-\gamma)\zeta+
    \frac{1}{2}\big(
    \vert {\bf H}[\zeta]\psi_1\vert^2-\gamma
    \vert \nabla \psi_1\vert^2\big)+
    {\mathcal N}(\zeta,\psi_1)=0
$$
where $\gamma=\rho_1/\rho_2$ and
\begin{eqnarray*}
    {\mathcal N}(\zeta,\psi_1)
    :=
    \frac{
    \gamma\big(G[\zeta]\psi_1+
    \nabla\zeta\cdot\nabla\psi_1\big)^2
    -\big(G[\zeta]\psi_1
    +\nabla\zeta\cdot {\bf H}[\zeta]\psi_1\big)^2
    }{2(1+\vert\nabla\zeta\vert^2)}.
\end{eqnarray*}
Taking the gradient of this equation and using (\ref{eq1}) then
gives the system of equations
\begin{equation}\label{eqdepart}
    \left\lbrace
    \begin{array}{l}
    \dt \zeta-G[\zeta]\psi_1=0,\vspace{1mm}\\
    \dt({\bf H}[\zeta]\psi_1-\gamma\nabla\psi_1)
    +g(1-\gamma)\nabla
    \zeta\\
    \indent\indent+
    \frac{1}{2}\nabla\big(
    \vert {\bf H}[\zeta]\psi_1\vert^2
    -\gamma\vert \nabla \psi_1\vert^2\big)+
    \nabla{\mathcal N}(\zeta,\psi_1)=0,
    \end{array}\right.
\end{equation}
for $\zeta$ and $\psi_1$.  This is the system of equations that will
be used in the next sections to derive asymptotic models.

\begin{remark}
    \label{rem1}
    More precise definitions of the
    operators $G[\zeta]$ and ${\bf H}[\zeta]$ will be presented
    in Subsection \ref{sectnondim} and in Section \ref{sectas}.

\end{remark}
\begin{remark}
    Setting $\rho_1=0$, and thus $\gamma=0$, in the above equations,
    one recovers the usual surface water-wave equations written in terms
    of $\zeta$ and $\psi$ as in \cite{CSS2,CSS,Zakharov}.
\end{remark}

\subsection{Non-Dimensionalization of the Equations}\label{sectnondim}

The asymptotic behaviour of (\ref{eqdepart}) is more transparent
when these equations are written in dimensionless variables.
Denoting by $a$ a typical amplitude of the deformation of the
interface in question, and by $\lambda$ a typical wavelenth, the
following dimensionless indendent variables
$$
    \widetilde{X}:=\frac{X}{\lambda},\quad
    \widetilde{z}:=\frac{z}{d_1},\quad
    \widetilde{t}:=\frac{t}{\lambda/\sqrt{gd_1},\quad},
$$
are introduced.  Likewise, we define the dimensionless unknowns
$$
    \widetilde{\zeta}:=\frac{\zeta}{a},\quad
    \widetilde{\psi}_1:=\frac{\psi_1}{a\lambda\sqrt{g/d_1}},
$$
as well as the dimensionless parameter`s
$$
    \gamma:=\frac{\rho_1}{\rho_2},\quad
    \delta:=\frac{d_1}{d_2},\quad
    \eps:=\frac{a}{d_1},\quad
    \mu:=\frac{d_1^2}{\lambda^2}.
$$
Though they are redundant, it is also notationally convenient to
introduce two other parameter`s $\eps_2$ and $\mu_2$ defined as
$$
    \eps_2=\frac{a}{d_2}=\eps\delta,\qquad
    \mu_2=\frac{d_2^2}{\lambda^2}=\frac{\mu}{\delta^2}.
$$
\begin{remark}
    The parameters $\eps_2$ and $\mu_2$ correspond to $\eps$ and
    $\mu$ with $d_2$ rather than $d_1$
    taken as the unit of length in the vertical direction.
\end{remark}
Before writing (\ref{eqdepart}) in dimensionless variables, a
dimensionless Dirichlet-Neumann operator $G^\mu[\eps\zeta]\cdot$ is
needed, associated to the non-dimensionalized upper fluid domain
$$
    \Omega_1=\{(X,z)\in \R^{d+1},-1+\eps{\zeta}(X)<z<0\}.
$$
Throughout the discussion, it will be presumed that this domain
remains connected, so there is a positive value $H_1$ such that
\begin{equation}\label{flotte1}
    1-\eps{\zeta}\geq H_1\quad\mbox{ on }\quad \R^d.
\end{equation}
\begin{definition}
    \label{def1}
    Let ${\zeta}\in W^{2,\infty}(\R^d)$ be such that
    (\ref{flotte1}) is satisfied  and let
    $\psi_1\in H^{3/2}(\R^d)$.
    If  $\Phi_1$ is the unique solution in
    $H^2(\Omega_1)$ of the boundary-value problem
    \begin{equation}
    \label{eqDN}
    \left\lbrace
    \begin{array}{l}
    \mu\Delta \Phi_1+\dz^2 \Phi_1=0\quad\mbox{ in }\Omega_1, \\
    \dz \Phi_1\,_{\vert_{z=0}}=0,
    \qquad \Phi_1\,_{\vert_{z=-1+\eps{\zeta}(X)}}=\psi_1,
    \end{array}\right.
    \end{equation}
 then $G^\mu[\eps{\zeta}]\psi_1\in H^{1/2}(\R^d)$ is defined by
    $$
    G^\mu[\eps{\zeta}]\psi_1=
    -\mu\eps\nabla{\zeta}\cdot \nabla
    \Phi_1\,_{\vert_{z=-1+\eps{\zeta}}}
    +\dz \Phi_1\,_{\vert_{z=-1+\eps{\zeta}}}.
    $$
\end{definition}
\begin{remark}
    Another way to approach $G^\mu$ is to define
     $$G^\mu[\eps\zeta]\psi_1=\sqrt{1+\eps^2\vert\nabla\zeta\vert^2}\partial_n\Phi_1\,_{\vert_{z=-1+\eps{\zeta}}}$$
     where
    $\partial_n\Phi_1\,_{\vert_{z=-1+\eps{\zeta}}}$ stands for the
    upper conormal derivative associated to the elliptic
    operator $\mu\Delta\Phi_1+\dz^2\Phi_1$.
\end{remark}
In the same vein, one may define a dimensionless operator
    ${\bf H}^{\mu,\delta}[\eps\zeta]\cdot$
associated to the non-dimensionalized lower fluid domain
$$
    \Omega_2=
    \{(X,z)\in \R^{d+1},-1-1/\delta<z<-1+\eps{\zeta}(X)\},
$$
where it is  assumed as in (\ref{flotte1})
that there is an $H_2  >
0$ such that
\begin{equation}\label{flotte2}
    1+\eps\delta {\zeta}\geq H_2\quad\mbox{ on }\quad \R^d.
\end{equation}
\begin{definition}
    \label{def2}
    Let ${\zeta}\in W^{2,\infty}(\R^d)$ be such that
    (\ref{flotte1}) and (\ref{flotte2}) are satisfied,
    and   suppose that  $\psi_1\in H^{3/2}(\R^d)$ is given.  If the function
    $\Phi_2$ is the unique solution (up to a constant) of the
    boundary-value problem
    \begin{equation}
    \label{eqjump}
    \left\lbrace
    \begin{array}{l}
    \mu\Delta \Phi_2+\dz^2 \Phi_2=0\quad\mbox{ in }\Omega_2, \\
    \dz \Phi_2\,_{\vert_{z=-1-1/\delta}}=0,
    \qquad \partial_n \Phi_2\,_{\vert_{z=-1+\eps{\zeta}(X)}}=
    \frac{1}{(1+\eps^2\vert\nabla\zeta\vert^2)^{1/2}}
    G^\mu[\eps\zeta]\psi_1,
    \end{array}\right.
    \end{equation}
    then the operator ${\bf H}^{\mu,\delta}[\eps\zeta]\cdot$ is defined on $\psi_1$ by
    $$
        {\bf H}^{\mu,\delta}[\eps\zeta]\psi_1
        =\nabla(\Phi_2\,_{\vert_{z=-1+\eps\zeta}})  \in H^{1/2}(\R^d).
    $$
\end{definition}
\begin{remark}
    In the statement above,
    $\partial_n \Phi_2\,_{\vert_{z=-1+\eps{\zeta}}}$ stands
    here for the upwards conormal derivative associated
    to the elliptic operator $\mu\Delta \Phi_2+\dz^2 \Phi_2$,
    $$
    \sqrt{1+\eps^2\vert\nabla\zeta\vert^2}
    \partial_n \Phi_2\,_{\vert_{z=-1+\eps{\zeta}}}
    =-\mu\eps \nabla\zeta\cdot\nabla \Phi_2\,_{\vert_{z=-1+\eps{\zeta}}}
    +\partial_z \Phi_2\,_{\vert_{z=-1+\eps{\zeta}}}.
    $$
    The Neumann boundary condition of (\ref{eqjump}) at the interface
    can also be stated as
    $\partial_n \Phi_2\,_{\vert_{z=-1+\eps{\zeta}}}
    =\partial_n \Phi_1\,_{\vert_{z=-1+\eps{\zeta}}}$.
\end{remark}
\begin{remark}\label{rmres}
    Of course, the solvability of (\ref{eqjump})
    requires the condition
    $\int_\Gamma \partial_n\Phi_2 d\Gamma=0$ (where $d\Gamma=\sqrt{1+\eps^2\vert\nabla\zeta\vert^2}dX$
    is the Lebesgue measure on  the surface $\Gamma=\{z=-1+\eps\zeta\}$).
    This is automatically satisfied thanks to the definition of
    $G^\mu[\eps\zeta]\psi_1$. Indeed, applying Green's identity to (\ref{eqDN}), one
    obtains
    $$
    \int_{\Gamma}  \partial_n\Phi_2d\Gamma
    =\int_{\Gamma}\partial_n\Phi_1d\Gamma
    =-\int_{\Omega_1} (\mu\Delta\Phi_1+\dz^2\Phi_1)=0.
    $$
\end{remark}
\begin{example}\label{ex1}
    The operators
    $G^\mu[\eps\zeta]\cdot$ and ${\bf
    H}^{\mu,\delta}[\eps\zeta]\cdot$ have explicit expressions
    when the interface is flat (i.e. when $\zeta=0$). In that case,
    taking the horizontal Fourier transform of the
    Laplace equations (\ref{eqDN}) and (\ref{eqjump})
    transforms them into ordinary differential equations with respect to $z$ which can
    easily be solved to obtain
    $$
    G^\mu[0]\psi=-\sqrt{\mu}\vert D\vert \tanh(\sqrt{\mu} \vert D\vert )\psi\quad\mbox{ and }\quad
    {\bf H}^{\mu,\delta}[0]\psi=-
    \frac{\tanh(\sqrt{\mu}\vert D\vert)}
    {\tanh(\frac{\sqrt{\mu}}{\delta}\vert D\vert)}\nabla\psi.
    $$
\end{example}
The equations (\ref{eqdepart}) can therefore be written in
dimensionless variables as
\begin{equation}
    \label{eqdepndim}
    \left\lbrace
    \begin{array}{lcl}
    \dsp \partial_{\widetilde{t}} \widetilde{\zeta} -
    \frac{1}{\mu}G^\mu[\eps\widetilde{\zeta}]\widetilde{\psi}_1&=&0,\\
    \dsp \partial_{\widetilde{t}} \big(
    {\bf H}^{\mu,\delta}[\eps\widetilde{\zeta}]\widetilde{\psi}_1-\gamma
    \nabla\widetilde{\psi}_1\big)
    +(1-\gamma)\nabla\widetilde{\zeta}& &\\
    \dsp \indent+
    \frac{\eps}{2}\nabla\big(\vert {\bf H}^{\mu,\delta}[{\eps\widetilde{\zeta}}]\widetilde{\psi}_1\vert^2-\gamma\vert\nabla\widetilde{\psi}_1\vert^2\big)+\eps\nabla {\mathcal N}^{\mu,\delta}(\eps\widetilde{\zeta},\widetilde{\psi}_1)&=&0,
    \end{array}\right.
\end{equation}
where ${\mathcal N}^{\mu,\delta}$ is defined for all pairs
$({\zeta},{\psi})$ smooth enough by the formula
\begin{eqnarray*}
    {\mathcal N}^{\mu,\delta}({\zeta},{\psi}):=
    \mu\frac{
    \gamma\big(\frac{1}{\mu}G^\mu[\zeta]\psi+
    \nabla\zeta\cdot\nabla\psi\big)^2
    -\big(\frac{1}{\mu}G^\mu[\zeta]\psi
    +\nabla\zeta\cdot{\bf H}^{\mu,\delta}[\zeta]\psi\big)^2
    }{2(1+\mu\vert\nabla\zeta\vert^2)}.
\end{eqnarray*}

Our work centers around the study of the asymptotics of the
non-dimensionalized equations (\ref{eqdepndim}) in various physical
regimes corresponding to different relationships among the
dimensionless parameter`s $\eps$, $\mu$ and $\delta$.
\begin{notation}
    The tildes which indicate the non-dimensional quantities will be systematically
    dropped henceforth.
\end{notation}
\begin{remark}
Linearizing the equations (\ref{eqdepndim}) around the rest state,
one finds the equations
$$
    \left\lbrace
    \begin{array}{lcl}
    \dsp \partial_{t} {\zeta}
    - \frac{1}{\mu}G^\mu[0]{\psi}_1&=&0,\\
    \dsp \partial_{{t}} \big(
    {\bf H}^{\mu,\delta}[0]{\psi}_1-\gamma\nabla{\psi}_1\big)
    +(1-\gamma)\nabla{\zeta}&= &0.
    \end{array}\right.
$$
The explicit formulas in Example \ref{ex1} thus allow one to
calculate the linearized dispersion relation
\begin{equation}\label{disprel}
    \omega^2=(1-\gamma)\frac{\vert{\bf k}\vert}{\sqrt{\mu}}\frac{\tanh(\sqrt{\mu}\vert {\bf k}\vert)\tanh(\frac{\sqrt{\mu}}{\delta}\vert{\bf k}\vert)}{\tanh(\sqrt{\mu}\vert {\bf k}\vert)+\gamma\tanh(\frac{\sqrt{\mu}}{\delta}\vert{\bf k}\vert)};
\end{equation}
corresponding to plane-wave solutions $e^{i{\bf k}\cdot X-i\omega
t}$.  In particular, the expected instability is found when
$\gamma>1$, corresponding to the case wherein the heavier fluid lies
over the lighter one.  One also checks that the classical dispersion
relation
$$\omega^2=\frac{1}{\sqrt{\mu}}\vert {\bf
k}\vert\tanh(\sqrt{\mu}\vert {\bf k}\vert)$$
for surface water waves
is recovered when $\gamma=0$ and $\delta=1$.
\end{remark}

\subsection{Principal Results}

The overall goal here is to propose model systems of equations for
the internal waves by obtaining the asymptotic form of the equations
(\ref{eqdepndim}) in various regimes corresponding to different
values of the parameters $\eps$, $\delta$ and $\mu$. All these
asymptotic models are ($1+d$)-dimensional systems coupling the
surface elevation $\zeta$ to the variable ${\bf v}$ defined to be
\begin{equation}\label{defv}
    {\bf v}:={\bf H}^{\mu,\delta}[\eps{\zeta}]{\psi}_1
    -\gamma
    \nabla{\psi}_1.
\end{equation}
(For the surface water-wave problem formally recovered by taking
$\gamma=0$ and $\delta=1$, ${\bf v}$ is the horizontal velocity
evaluated at the free surface).  We will often refer to ${\bf v}$ as
the velocity variable, though its precise interpretation will vary.
Note that ${\bf v}$ is essentially the gradient of the second canonical 
variable 
in the hamiltonian formulation of (\ref{eqdepndim}), (see for instance \cite{BB}). 

It will be rigorously established that the internal-wave equations
(\ref{eqdepndim}) are \emph{consistent} with the asymptotic models
for $(\zeta,{\bf v})$ derived in this paper in the following precise
sense.
\begin{definition}\label{defcons}
    The internal wave equations (\ref{eqdepndim}) are \emph{consistent}
    with a system $S$
    of $d+1$  equations for
    $\zeta$ and ${\bf v}$ if for all sufficiently smooth
    solutions $(\zeta,\psi_1)$
    of (\ref{eqdepndim})
such that
    (\ref{flotte1}) and (\ref{flotte2}) are satisfied,
    the pair $(\zeta,{\bf v}={\bf H}^{\mu,\delta}[\eps{\zeta}]{\psi}_1
    -\gamma
    \nabla{\psi}_1)$ solves $S$ up to a small residual called
    the {\it precision} of the asymptotic model.
\end{definition}

\begin{remark}
    It is worth emphasis that above definition does not require the well-posedness of
    the internal wave equations (\ref{eqdepndim}).  Indeed, these
    can be subject to Kelvin-Helmholtz type instabilities (see for instance
    \cite{BHL} and \cite{ITT}), although one might expect a ``stability of the
    instability" result even in the face of such instabilities (see \cite{PC}).
    {\it Consistency} is only concerned with the
    properties of smooth solutions to the system
    (which do exist in the classical configuration of the Kelvin-Helmholtz
    problem, even when instabilities manifest themselves; see e.g. \cite{SSBF,SS}).
    In fact, the two-layer water-wave system is
    known to be well-posed in
    Sobolev spaces in the presence of surface tension \cite{ITT}. In consequence, one
    could simply add a small amount of surface tension at the interface between the
two homogeneous layers to put oneself in a well-posed situation. The
resulting analysis would be exactly the same and would, in fact,
lead to the same asymptotic models. (Such an approach is used in
\cite{OI} for
    the Benjamin-Ono equation).  As the resulting model systems do
    not change, such a regularization has been eschewed here.
\end{remark}

Here is a summary of the different asymptotic regimes investigated
in this paper. It is convenient to organize the discussion around
the parameters $\eps$ and $\eps_2=\eps\delta$ (the nonlinearity, or
amplitude, parameters for the upper and lower fluids, respectively),
and in terms of $\mu$ and $\mu_2=\frac{\mu}{\delta^2}$ (the
long-wavelength parameters for the upper and lower fluids).  Notice
that the assumptions made about $\delta$ are therefore implicit.

The interfacial wave is said to be of \emph{small amplitude} for the
upper fluid layer (resp. the lower layer) if $\eps\ll 1$ (resp.
$\eps_2\ll1$) and the upper (resp. lower) layer is said to be
\emph{shallow} if $\mu\ll 1$ (resp. $\mu_2\ll1$). This terminology
is consistent with the usual one for surface water waves (recovered
by taking $\rho_1=0$ and $\delta=1$). In the discussion below, the
notation \emph{regime 1/regime 2} means that the wave motion is such
that the upper layer is in regime 1 (small amplitude or shallow
water) and the lower one is in regime 2.
\begin{enumerate}
    \item The small-amplitude/small-amplitude regime: $\eps\ll1$, $\eps_2\ll1$.
    This regime corresponds to interfacial deformations which are small
for both the upper and lower fluid domains. Various sub-regimes are
defined  by making further assumptions about the size of $\mu$ and
$\mu_2$.
    \begin{enumerate}
    \item The Full Dispersion /Full Dispersion (FD/FD) regime:
    $\eps\sim\eps_2\ll1$ and $\mu\sim\mu_2=O(1)$
    (and thus $\delta\sim 1$).
    In this regime, investigated in
    \S \ref{ILW-ILW}, the shallowness parameters are not small for
    either of the fluid domains, and the full dispersive effects
    must therefore be kept for both regions; the asymptotic model
    corresponding to this situation
    is given in (\ref{eqFDFD}).
    \item The Boussinesq / Full dispersion (B/FD) regime:
    $\mu\sim\eps\ll 1$,
    $\mu_2\sim 1$. This regime is studied in \S \ref{B-ILW} and
    corresponds to the case where the flow has a Boussinesq structure
    in the upper part (and thus dispersive effects of the same
    order as nonlinear effects), but with a shallowness parameter
    not small in the lower fluid domain. This configuration occurs
    when $\delta^2\sim \eps$, that is, when the lower region is much
    larger than the upper one. A further analysis of the asymptotic
    model yields a three-parameter family of equivalent systems
    (see (\ref{eqB-FD}) below).
    \item The Boussinesq/Boussinesq (B/B) regime:
    $\mu\sim\mu_2\sim\eps\sim\eps_2\ll1$.
    In this regime, investigated in \S \ref{B-B}, one has $\delta\sim 1$
    and the flow has a Boussinesq structure in
    both the upper and lower fluid domains. Here again, a three-parameter
    family of asymptotic systems is obtained (see (\ref{eqB-B}) below).
    \end{enumerate}

    \item The Shallow Water/Shallow Water (SW/SW) regime:
    $\mu\sim \mu_2\ll1$.
    This regime, which allows relatively large interfacial amplitudes
    ($\eps\sim\eps_2=O(1)$), does not belong to the
    regimes singled out above. The structure of the flow is
    then of shallow water type in both regions; in particular,
    the asymptotic model (see \S \ref{SW-SW}) is a nonlinear,
    but non-dispersive system, given in (\ref{eqSWSW}),
    which degenerates into the
    usual shallow water equations when $\gamma=0$ and $\delta=1$.
    It is very interesting in this case that a non-local term arises
    when $d=2$.  Such a nonlocal
    term does not appear in the one-dimensional case, nor in the
    two-dimensional shallow water equations for surface waves.
    \item The Shallow Water/Small Amplitude (SW/SA) regime:  $\mu\ll1$
    and $\eps_2\ll 1$. In this regime, the upper layer is
    shallow (but with possibly large surface deformations), and
    the surface deformations are small for the lower layer (but it
    can be deep). Various sub-regimes arise in this case also.
    \begin{enumerate}
    \item The Shallow Water/Full dispersion (SW/FD) regime:
    $\mu\sim \eps_2^2\ll 1$, $\eps\sim\mu_2\sim 1$. This
    regime is investigated in \S \ref{SW-FD}. The dispersive
    effects are negligible  in the upper fluid, but the
    full dispersive effects must be
    kept in the lower one (see system (\ref{eqSWFD}) below).
    \item The Intermediate Long Waves (ILW) regime:
     $\mu\sim \eps^2\sim\eps_2\ll1$, $\mu_2\sim 1$. In this regime,
    the interfacial deformations are also small for the upper fluid
    (which is not the case in the SW/FD regime).  This allows some
    simplifications, as shown in \S \ref{ILW}. It is also
    possible (see (\ref{eqILW}))
    to derive a one-parameter family of equivalent systems.
    \item The Benjamin-Ono (BO) regime:
    $\mu\sim \eps^2\ll1$, $\mu_2=\infty$. A formal study of this
    regime is performed in \S \ref{BO}.  It is shown in particular
    how to recover the Benjamin-Ono equation as the unidirectional
    limit in the one-dimensional case $d=1$.  The Benjamin-Ono
    equation is also shown to be a particular case of a one-parameter
    family of \emph{regularized Benjamin-Ono equations}, given in
    (\ref{eqBO}).
    \end{enumerate}
\end{enumerate}

The range of validity of these regimes is summarized in the
following table.

\begin{center}
\begin{tabular}{|*{3}{l|}}
\hline
     & $\eps=O(1)$ & $\eps\ll1$   \\
    \hline
    $\mu=O(1)$ & Full equations & $\delta\sim 1$: FD/FD eq'ns\\
    \hline
    $\mu\ll1$ & $\delta\sim 1$: SW/SW eq'ns & $\mu\sim\eps$ and $\delta^2\sim\eps$: B/FD eq'ns\\
    &  $\delta^2\sim \mu\sim\eps_2^2$: SW/FD eq'ns & $\mu\sim\eps$ and $\delta\sim1$: B/B eq'ns\\
    &  &$\delta^2\sim \mu\sim\eps^2$: ILW eq'ns \\
    &  &$\delta=0$ and $\mu\sim\eps^2$: BO eq'ns  \\
    \hline
    \end{tabular}
\end{center}

\begin{remark}
    The small amplitude/shallow water regime is not investigated
    here.  It corresponds to the situation where the upper fluid domain
    is much larger than the lower one, which is more of an
    atmospheric configuration than an oceanographic case.
\end{remark}

\section{Asymptotic Expansions of the Operators}\label{sectas}

In this section, asymptotic expansions are given of the central
operators defined in the Introduction.  The discussion begins with
the Dirichlet-Neumann operator.

\subsection{Asymptotic expansion of the Dirichlet-Neumann Operator
$G^\mu[\eps\zeta]\cdot$}

The following lemma connects $\zeta$ with the vertically integrated
horizontal velocity via the Dirichlet-Neumann operator
$G^\mu[\eps\zeta]\cdot$.
\begin{lemma}\label{lmvert}
    Let $\zeta\in W^{2,\infty}(\R^d)$ be such that (\ref{flotte1})
    is satisfied and let $\psi \in H^{3/2}(\R^d)$
    and $\Phi_1$ be the solution
    of (\ref{eqDN}) with $\psi_1=\psi$. If $V^\mu$ is defined by
    $$V^{\mu}[\eps\zeta]\psi:=\int_{-1+\eps\zeta}^0
    (\sqrt{\mu}\nabla\Phi_1)dz, $$
    then it follows that
    $$
    G^\mu[\eps\zeta]\psi =
    \sqrt{\mu}\nabla\cdot (V^{\mu}[\eps\zeta]\psi).
    $$
\end{lemma}
\begin{proof}
Let $\varphi\in C^\infty_c(\R^d)$ be a test function.  Using Green's
identity, and with the notation of Remark \ref{rmres}, one obtains
\begin{eqnarray*}
    \int_{\R^d}G^\mu[\eps\zeta]\psi \varphi&=&
    \int_{\Gamma}\partial_n\Phi_1\varphi d\Gamma\\
    &=&-\int_{\Omega_1}(\sqrt{\mu}\nabla)\Phi_1\cdot (\sqrt{\mu}\nabla)\varphi\\
    &=&-\int_{\R^d}\int_{-1+\eps\zeta}^0(\sqrt{\mu}\nabla\Phi_1)dz\cdot \sqrt{\mu}\nabla \varphi.
\end{eqnarray*}
Defining $V^{\mu}[\eps\zeta]\psi$ as in the statement of the lemma,
it transpires that
$$
    \int_{\R^d}G^{\mu}[\eps\zeta]\varphi
    =-\sqrt{\mu}\int_{\R^d}V^{\mu}[\eps\zeta]\psi\cdot \nabla\varphi
    =\sqrt{\mu}\int_{\R^d}\nabla\cdot (V^{\mu}[\eps\zeta]\psi)\varphi.
$$
Since the above identity is true for all $\varphi\in C^\infty_c(\R^d)$,
the result follows.
\end{proof}
\begin{remark}
    In \S \ref{sectdassmall} and \S \ref{sectlw} below,
    asymptotic expansions are obtained of $V^{\mu}[\eps\zeta]\psi$ in terms of $\eps$
    and $\mu$, respectively. Because of Lemma \ref{lmvert}, asymptotic expansions of
    $G^\mu[\eps\zeta]\psi $ then follow immediately.
\end{remark}
\subsubsection{Asymptotic Expansion of
$V^\mu[\eps\zeta]\cdot$ when
${\mathbf\eps}\ll1$}\label{sectdassmall}

When $\eps\ll1$, the approach to obtaining an asymptotic expansion
of $V^\mu[\eps\zeta]\psi$ is to make a Taylor expansion in terms of
the interface deformation around the rest state, {\it viz.}
$$
    V^\mu[\eps\zeta]\psi =V^\mu[0]\psi
    +\eps (d_0(V^\mu[\cdot])\zeta)\psi
    +\cdots.
$$
(Note, however, that the expansion of $V^\mu[\eps\zeta]\psi$ itself,
and not only the consequent expansion of $G^\mu[\eps\zeta]\psi$, is
needed so that the elliptic estimate of Proposition \ref{expgal} can
be used in the proof of Corollary \ref{coro2}).
\begin{proposition} \label{prop2}
    Let $s>d/2$ and $\zeta\in H^{s+3/2}(\R^d)$ be such that
    (\ref{flotte1}) is satisfied. Then for
    $\psi$ such that $\nabla\psi\in H^{s+1/2}(\R^d)$,
    the inequality
    \begin{eqnarray*}
    \lefteqn{\big\vert V^{\mu}[\eps\zeta]\psi-\big[
    {\mathcal T}_{0,\mu}\nabla \psi
    +\eps \sqrt{\mu} (-\zeta+{\mathcal T}_{1,\mu}[\zeta])
    \nabla\psi\big]\big\vert_{H^s}}\\
    & &\indent\indent\leq \eps^2
    C(\frac{1}{H_1},\eps\sqrt{\mu},\vert\zeta\vert_{H^{s+3/2}}
    ,\vert\nabla\psi\vert_{H^{s+1/2}}),
    \end{eqnarray*}
    holds for all
    $\eps\in [0,1]$ and $\mu>0$, where
    ${\mathcal T}_{0,\mu}=\frac{\tanh(\sqrt{\mu}\vert D\vert)}{\vert
D\vert}$, ${\mathcal T}_{1,\mu}[\zeta]=
    -\nabla{\mathcal T}_{0,\mu} (\zeta {\mathcal T}_{0,\mu}\nabla^T)
    $, and  $V^{\mu}[\eps\zeta]\psi$ is as defined in
    Lemma \ref{lmvert}
    (so that $G^{\mu}[\eps\zeta]\psi=\sqrt{\mu}\nabla\cdot V^{\mu}[\eps\zeta]\psi$).
\end{proposition}
The key point in the proof is an \emph{explicit} formula of the
derivative of the mapping $\zeta\mapsto V^\mu[\eps\zeta]\psi$, which
generalizes the formula obtained in \cite{Lannes} for the shape
derivative of Dirichlet-Neumann operators.  This interesting
technical point is the subject of the next lemma.
\begin{lemma}\label{lmformula}
    Let $s>d/2$ and suppose that $\psi$ is such that $\nabla\psi\in H^{s+1/2}(\R^d)$.
    The mapping $H^{s+3/2}(\R^d)\ni \zeta\mapsto V^\mu[\eps\zeta]\psi\in H^{s+1/2}(\R^d)^d$ is differentiable. Moreover, for all $\zeta,\zeta'\in H^s(\R^d)$,
    the derivative of $V^\mu[\eps\cdot]\psi$ at $\zeta$ in the
    direction $\zeta'$ is given by the formula
    $$
    d_{\zeta} (V^\mu[\eps\cdot]\psi)\zeta'=
    -\eps V^\mu[\eps\zeta](\zeta'Z^\mu[\eps\zeta]\psi)
    -\eps \zeta'
    \big(\sqrt{\mu}\nabla\psi
    -\eps\sqrt{\mu}\nabla\zeta Z^\mu[\eps\zeta]\psi\big),
    $$
    where $Z^\mu[\eps\zeta]\psi:=\frac{1}{1+\eps^2\mu\vert\zeta\vert^2}
    (G^\mu[\eps\zeta]\psi+\eps\mu\nabla\zeta\cdot\nabla\psi)$.
\end{lemma}
\begin{proof}[Proof of the Lemma]
First, define another Dirichlet-Neumann operator ${\mathcal
G}^\mu[\eps\zeta]\cdot$ by
\begin{equation}\label{DNbis}
    {\mathcal G}^\mu[\eps\zeta]\psi={\bf e_z}\cdot P^\mu[\eps\zeta]\nm\Phi_{\vert_{z=0}},
\end{equation}
where $\Phi$ solves
\begin{equation}\label{eqnP}
    \left\lbrace
    \begin{array}{l}
    \nm\cdot P^\mu[\eps\zeta]\nm\Phi=0 \quad\mbox{ in }\quad -1<z<0, \\
    \Phi\,_{\vert_{z=0}}=\psi,\qquad
    \dz \Phi\,_{\vert_{z=-1}}=0,
    \end{array}\right.
\end{equation}
and where
$$\displaystyle P^\mu[\eps\zeta]=\left(\begin{array}{cc}
(1+\eps\zeta) I_{d\times d} & \eps\sqrt{\mu}(z+1)\nabla\zeta\\
-\eps\sqrt{\mu}(z+1)\nabla\zeta^T &
\frac{1+\eps^2\mu\vert\nabla\zeta\vert^2}{1+\eps\zeta}\end{array}\right).$$
This operator is the classical Dirichlet-Neumann operator often used
for the study of the surface water-wave equations and for which an
explicit expression exists for the derivative of the mapping
$\zeta\mapsto {\mathcal G}^\mu[\eps\zeta]\psi$ (see, e.g. Theorem
3.20 of \cite{Lannes} and Theorem 3.1 of \cite{AL}). Studying the
transformation of the fluid domain into the flat strip $-1<z<0$
(flattening of the domain) reveals that
$G^\mu[\eps\zeta]\psi=-{\mathcal G}^\mu[-\eps\zeta]\psi$ (see
Proposition 2.7 of \cite{AL} and Section \ref{asH} below where the
same kind of transformation is performed). It will be convenient to
consider the operator ${\mathcal G}^\mu[-\eps\zeta]\cdot$ rather
than $G^\mu[\eps\zeta]\cdot$ because this allow us to take over
intact some elements of the proof of Theorem 3.20 in \cite{Lannes}.
Moreover, for the sake of clarity, we take $\eps=\mu=1$ in this
proof and leave to the reader the straighforward modifications for
the general case. The
proof is divided into 5 steps.\\
{\it Step 1.} One has that $\cG\psi=-\nabla\cdot(\cV)$, with $
    \cV=\int_{-1}^0 P_{I}[\zeta]\nabla_{X,z}\Phi dz,
$ and where $P_{I}[\zeta]$ is the $d\times (d+1)$ matrix obtained by
taking the last row off $P[\zeta]$. The proof of this result is more
or less identical
to the proof of Lemma \ref{lmvert}.\\
{\it Step 2.} Denoting by ${\mathcal V}'$ the derivative of
${\mathcal V}[\cdot]\psi$ at $\zeta$ and in the direction $\zeta'$,
one computes
$$
    {\mathcal V}'=
    \int_{-1}^0 P'_I\nabla_{X,z}\Phi dz+\int_{-1}^0
    P_{I}[\zeta]\nabla_{X,z}\Phi' dz,
$$
where $P'_I$ and $\Phi'$ stand, respectively, for the derivative at $\zeta$
and in the direction $\zeta'$ of the mappings
$\zeta\mapsto P_I[\zeta]$ and $\zeta\mapsto \Phi$.\\
{\it Step 3.} Defining $\displaystyle
\chi=(z+1)\frac{\zeta'}{1+\zeta}\dz\Phi$, one has
$$
    \int_{-1}^0P_I[\zeta]\nabla_{X,z}(\Phi'-\chi)
    =-\cV(\zeta'\cZ\psi),
$$
with $\displaystyle
\cZ\psi=\frac{\cG\psi+\nabla\zeta\cdot\nabla\psi}{1+\vert\nabla\zeta\vert^2}$.
To prove this result, first remark that $w:=\Phi'-\chi$ solves the
boundary-value problem
$$
    \left\lbrace
    \begin{array}{l}
    \nm\cdot P^\mu[\eps\zeta]\nm w=0 \quad\mbox{ in }\quad -1<z<0, \\
    w \,_{\vert_{z=0}}=-\zeta'\cZ\psi,\qquad
    \dz w\,_{\vert_{z=-1}}=0,
    \end{array}\right.
$$
as a consequence of Lemma 3.22 of \cite{Lannes}. The result then
follows directly from the definition of $\cV\cdot $.\\
{\it Step 4.} The identity
$$
    \int_{-1}^0(P'_I\nabla_{X,z}\Phi+P_I[\zeta]\nabla_{X,z}\chi)dz
    =\zeta'\big(\nabla\psi-\cZ\psi\nabla\zeta\big).
$$
also holds.  To establish this, first compute that
$$
    P'_I\nabla_{X,z}\Phi+P_I[\zeta]\nabla_{X,z}\chi=
    \zeta'\dz\big((z+1)\nabla\Phi\big)
    -\nabla\zeta\dz\Big(\frac{(z+1)^2}{h}\zeta'\dz\Phi\Big).
$$
The result then follows upon integrating with respect to $z$.\\
{\it Step 5.} It now remains simply to put together the pieces.  It
is deduced from Steps 2-4 that
$$
    {\mathcal V}'=\zeta'\big(\nabla\psi-\cZ\psi\nabla\zeta\big)-\cV(\zeta'\cZ\psi).
$$
The result then follows from the observation that if $V[\zeta]\psi$
is as defined in Lemma \ref{lmvert}, then
 $V[\zeta]\psi={\mathcal V}[-\zeta]\psi$.
\end{proof}
\begin{proof}[Proof of the Proposition]
A second order Taylor expansion reveals that
$$
    V^\mu[\eps\zeta]\psi=
    V^\mu[0]\psi+d_0(V^\mu[\eps\cdot]\psi)\zeta+
    \int_0^1 (1-z)d^2_{z\zeta}(V^\mu[\eps\cdot]\psi)(\zeta,\zeta)dz.
$$
Lemma \ref{lmformula} therefore implies that
\begin{eqnarray*}
    \lefteqn{V^\mu[\eps\zeta]\psi=
    V^\mu[0]\psi}\\& &-\eps V^\mu[0](\zeta G^\mu[0]\psi)
    -\eps\sqrt{\mu}\zeta \nabla\psi
    +\int_0^1
    (1-z)d^2_{z\zeta}(V^\mu[\eps\cdot]\psi)(\zeta,\zeta)dz.
\end{eqnarray*}

We saw in Example \ref{ex1} that $G^\mu[0]\psi=-\sqrt{\mu}\vert
D\vert \tanh(\sqrt{\mu} \vert D\vert )\psi$.  Similarly, one can
check that $V^\mu[0]\psi=\frac{\tanh(\sqrt{\mu}\vert D\vert)}{\vert
D\vert}\nabla \psi$. The proof of the proposition is now clear after
appreciating that
$$
    \big\vert \int_0^1
    (1-z)d^2_{z\zeta}(V^\mu[\eps\cdot]\psi)(\zeta,\zeta)dz\big\vert_{H^s}
    \leq
    \eps^2C(\frac{1}{H_1},\eps\sqrt{\mu},\vert\zeta\vert_{H^{s+3/2}}
    ,\vert\nabla\psi\vert_{H^{s+1/2}}),
$$
a fact which is obtained exactly as in Proposition 3.3 of \cite{AL}.
\end{proof}

\subsubsection{Asymptotic Expansion of $V^\mu[\eps\zeta]\cdot$ for Large-Amplitude
Waves and Shallow Depth ($\eps=O(1)$ and $\mu\ll 1$)}\label{sectlw}

For larger amplitude waves, the expansion of the Dirichlet-Neuman
operator $G^\mu[\eps\zeta]\psi$ (and also of $V^\mu[\eps\zeta]\psi$)
around the rest state no longer provides an accurate approximation.
 However, if $\mu\ll 1$, which is what we have earlier called
  the shallow water regime for the upper fluid,
 it is possible to obtain an expansion of $V^\mu[\eps\zeta]\psi$
(and thus of $G^\mu[\eps\zeta]\psi=\sqrt{\mu}\nabla\cdot
V^\mu[\eps\zeta]\psi$) with respect to $\mu$ which is uniform with
respect to $\eps\in[0,1]$.
\begin{proposition}
    \label{prop1}
    Let $s>d/2$ and $\zeta\in H^{s+3/2}(\R^d)$. Then for all
    $\mu\in (0,1)$ and
    $\psi$ such that $\nabla\psi\in H^{s+5/2}(\R^d)$, one has
    $$
    \big\vert \sqrt{\mu}V^{\mu}[\eps\zeta]\psi-
    \mu (1-\eps\zeta)\nabla\psi\big\vert_{H^{s}}\leq \mu^2
    C(\vert \zeta\vert_{H^{s+3/2}},\vert \nabla\psi\vert_{H^{s+5/2}}),
    $$
    uniformly with respect to
    $\eps\in [0,1]$), where $V^{\mu}[\eps\zeta]\psi$ is as defined in
    Lemma \ref{lmvert}
    (so that $G^{\mu}[\eps\zeta]\psi=\sqrt{\mu}\nabla\cdot V^{\mu}[\eps\zeta]\psi$).
\end{proposition}
\begin{remark}\label{remB}
    As in Prop. 3.8 of \cite{AL}, one can carry out the expansion explicitly to
    second order in $\mu$, thereby obtaining
    $$
    \sqrt{\mu}V^{\mu}[\eps\zeta]\psi=
    \mu (1-\eps\zeta)\nabla\psi+\frac{\mu^2}{3}\Delta\nabla\psi
    +O(\mu^3,\eps\mu^2).
    $$
\end{remark}
\begin{proof}
Recall that $G^\mu[\eps\zeta]\psi=-{\mathcal
G}^\mu[-\eps\zeta]\psi$, where ${\mathcal G}^\mu[\eps\zeta]\cdot$ is
defined in (\ref{DNbis}), and that ${\mathcal G}^\mu[\eps\zeta]\psi
=-\sqrt{\mu}\nabla\cdot V^\mu[-\eps\zeta]\psi$ (see the proof of
Lemma \ref{lmformula}). Proposition 3.8 of \cite{AL} shows that
$$
    \big\vert {\mathcal G}^{\mu}[\eps\zeta]\psi-
    \nabla\cdot(-\mu (1+\eps\zeta)\nabla\psi)\big\vert_{H^{s}}\leq \mu^2
    C(\vert \zeta\vert_{H^{s+3/2}},\vert
    \nabla\psi\vert_{H^{s+5/2}}).
$$
An obvious adaptation of the proof shows that the estimate given in
the statement of the Proposition can be obtained in the same way.
\end{proof}

\subsection{Asymptotic Expansions of
${\mathbf H^{\mu,\delta}}[\eps\zeta]\cdot$}\label{asH}

Attention is now turned to the interface operator
 ${\mathbf H^{\mu,\delta}}[\eps\zeta]\cdot$.\\
 The boundary-value problem (\ref{eqjump}) plays a key role in the analysis of the
operator ${\mathbf H^{\mu,\delta}}[\eps\zeta]\cdot$.  The analysis
of this problem is easier if we first transform it into a
variable-coefficient, boundary-value problem on the flat strip
${\mathcal S}:=\R^d\times (-1,0)$ using the diffeomorphism
$$
    \sigma: \begin{array}{ccc}
    {\mathcal S}& \to & \Omega_2\\
    (X,z)&\mapsto &
    \sigma(X,z):=(X,(1+\eps\delta)\frac{z}{\delta}+(-1+\eps\zeta)).
        \end{array}
$$
As shown in Proposition 2.7 of \cite{Lannes} (see also \S 2.2 of
\cite{AL}), $\Phi_2$ solves (\ref{eqjump}) if and only if
$\underline{\Phi}_2:=\Phi_2\circ \sigma$ solves
\begin{equation}\label{newbvp}
    \left\lbrace
    \begin{array}{l}
    \nabla^{\mu_2}_{X,z}\cdot Q^{\mu_2}[\eps_2\zeta]\nabla^{\mu_2}_{X,z}
    \underline{\Phi}_2=0 \qquad\mbox{ in }{\mathcal S},
    \vspace{1mm}\\
    \partial_n \underline{\Phi}_2\,_{\vert_{z=0}}=\frac{1}{\delta}
    G^\mu[\eps\zeta]\psi_1
    ,\qquad
    \partial_n \underline{\Phi}_2\,_{\vert_{z=-1}}=0,
    \end{array}\right.
\end{equation}
with
$$
    Q^{\mu_2}[\eps_2\zeta]=
    \left(\begin{array}{cc}
    (1+\eps_2\zeta)I_{d\times d} & -\sqrt{\mu_2}\eps_2 (z+1)\nabla\zeta\\
     -\sqrt{\mu_2}\eps_2 (z+1)\nabla\zeta^T &
    \frac{1+\mu_2\eps_2^2(z+1)^2\vert\nabla\zeta\vert^2}{1+\eps_2\zeta}
          \end{array}\right),
$$
and where, as before, $\eps_2=\eps\delta$,
$\mu_2=\frac{\mu}{\delta^2}$, and
$\nabla_{X,z}^{\mu_2}=(\sqrt{\mu_2}\nabla,\dz)^T$.
\begin{remark}
    As always in the present exposition, $\partial_n\underline{\Phi}_2$
    stands for the \emph{upward} conormal derivative associated
    to the elliptic operator involved in the boundary-value problem,
    $$
    \partial_n \underline{\Phi}_2\,_{\vert_{z=0\mbox{ \small{or} }z=-1}}={\bf e_z}\cdot Q^{\mu_2}[\eps_2\zeta]\nabla_{X,z}^{\mu_2}\underline{\Phi}_2\,_{\vert_{z=0\mbox{ \small{or} }z=-1}},
    $$
    where ${\bf e_z}$ is the upward-pointing unit vector along the vertical axis.
\end{remark}
An asymptotic expansion of
\begin{equation}\label{eqdefH}
    {\mathbf H^{\mu,\delta}}[\eps\zeta]\psi_1
    =\nabla(\underline{\Phi}_2\,_{\vert_{z=0}}),
\end{equation}
is obtained by finding an approximation $\underline{\Phi}_{app}$ to
the solution of (\ref{newbvp}) and then using the formal
relationship ${\mathbf H^{\mu,\delta}}[\eps\zeta]\psi_1\sim \nabla
(\underline{\Phi}_{app} \,_{\vert_{z=0}})$. This procedure is
justified in the following proposition, whose proof is postponed to
Appendix \ref{appA} so as not to interrupt the flow of the
development.  The proposition is used in both \S \ref{smalleps} and
\S \ref{largeeps} to give explicit asymptotic expansions of
${\mathbf H^{\mu,\delta}}[\eps\zeta]\psi_1$.  To state the result,
it is useful to have in place the spaces
$$
    H^{s,k}({\mathcal S})=\{f\in {\mathcal D}'(\overline{{\mathcal
    S}}): \Vert f\Vert_{H^{s,k}}<\infty\}
$$
for $s \in \R$ and $k \in \N$, where $\Vert
f\Vert_{H^{s,k}}=\sum_{j=0}^k\Vert \Lambda^{s-j}\dz^j f\Vert$.
\begin{proposition}\label{expgal}
    Let $s_0>d/2$, $s\geq s_0+1/2$, and $\zeta\in H^{s+3/2}(\R^d)$
    be such that
    (\ref{flotte1}) and (\ref{flotte2}) are satisfied (the interface
    does not touch the horizontal boundaries).
    If ${\bf h}\in H^{s+1/2,1}({\mathcal S})^{d+1}$
    and $V\in H^{s+1}(\R^d)^d$ are given, then the
    boundary-value problem
    \begin{equation}\label{eqlm}
    \left\lbrace
    \begin{array}{l}
    \nabla^{\mu_2}_{X,z}\cdot Q^{\mu_2}[\eps_2\zeta]\nabla^{\mu_2}_{X,z}
    u=\nabla_{X,z}^{\mu_2}\cdot {\bf h} \qquad\mbox{ in }{\mathcal S},
    \vspace{1mm}\\
    \partial_n u_{\vert_{z=0}}=\sqrt{\mu_2}\nabla\cdot V+
    {\bf e_z}\cdot {\bf h}_{\vert_{z=0}}
    ,\qquad
    \partial_n u_{\vert_{z=-1}}={\bf e_z}\cdot {\bf h}_{\vert_{z=-1}}
    \end{array}\right.
    \end{equation}
    admits a unique solution $u$.  Moreover, the solution $u$ obeys the
    inequality
    $$
    \big\vert
    \nabla u_{\vert_{z=0}}\big\vert_{H^s}\leq
    \frac{1}{\sqrt{\mu_2}}
    C(\frac{1}{H_2},\eps_2^{max},\mu_2^{max},\vert\zeta\vert_{H^{s+3/2}})
    \big(\Vert {\bf h}\Vert_{H^{s+1/2,1}}+\vert V\vert_{H^{s+1}}\big),
    $$
    uniformly with respect to $\eps_2\in [0,\eps_2^{max}]$ and
    $\mu_2\in (0,\mu_2^{max})$.
\end{proposition}
\begin{remark}
     In the case of a flat interface ($\zeta=0$),
    Example \ref{ex1} shows that $\frac{1}{\delta}G^{\mu}[0]\psi_1=
    \sqrt{\mu_2}\nabla\cdot V$ with
    $V=\frac{\nabla}{\vert D\vert}\tanh(\sqrt{\mu}\vert D\vert)\psi_1$.
    Consequently, (\ref{newbvp}), (\ref{eqdefH}) and
    Proposition \ref{expgal} (with ${\bf h}=0$) show that
    $$
    \big\vert {\bf H}^{\mu,\delta}[0]\psi_1\big\vert_{H^s}
    \lesssim \big\vert
    \frac{\tanh(\sqrt{\mu}\vert D\vert)}{\sqrt{\mu_2}\vert D\vert}
    \nabla\psi_1\big\vert_{H^{s+1}}\lesssim
    \delta\vert \nabla\psi_1\vert_{H^{s+1}},
    $$
    which is exactly the estimate one could have deduced from
    the explicit expression for ${\bf H}^{\mu,\delta}[0]\cdot$
    given in Example \ref{ex1} (except that using the latter approach gives
    an estimate
    in $H^s$ rather than in $H^{s+1}$.  The $H^s$-type result does
    not in fact carry over to
     the general case of non-flat interfaces).
\end{remark}
\begin{remark}
     Suppose we take ${\bf h}=0$ and $V=V^\mu[\eps\zeta]\psi$
    in Proposition 3.  By Lemma \ref{lmvert}, one has
    $G^\mu[\eps\zeta]\psi=\sqrt{\mu}\nabla\cdot V^\mu[\eps\zeta]\psi$,
    and so
    it follows that
    $\nabla u_{\vert_{z=0}}={\bf H}^{\mu,\delta}[\eps\zeta]\psi$.
    The Proposition thus provides an estimate of
    the operator norm of ${\bf H}^{\mu,\delta}[\eps\zeta]$.
\end{remark}

\subsubsection{The Small-Amplitude/Small-Amplitude Regime: $\eps\ll 1$, $\eps_2\ll 1$}\label{smalleps}

In this regime, it is assumed that the interface deformations are of
small amplitude for both the upper and lower fluids. The asymptotic
expansion of the operator $H^{\mu,\delta}[\eps\zeta]$ is thus made
in terms of $\eps$ and $\eps_2=\eps\delta$.  We proceed by first
constructing formally an approximate solution
$\underline{\Phi}_{app}$ to (\ref{newbvp}) in the form
$$
    \underline{\Phi}_{app}={\Phi^{(0)}}+\eps_2{\Phi^{(1)}}.
$$
This formal approximation is then justified rigorously in Corollary
\ref{coro2} below.

We may write from the expression for $Q^{\mu_2}[\eps_2\zeta]$, 
$$
    \nm\cdot Q^{\mu_2}[\eps_2\zeta]\nm=\dm+\eps_2 \nm\cdot Q_1\nm
    +\eps_2^2\nm\cdot Q_2\nm,
$$
with
$$
    Q_1=\left(
    \begin{array}{cc} \zeta I_{d\times d} & -\sqrt{\mu_2}(z+1)\nabla\zeta \\
    -\sqrt{\mu_2}(z+1)\nabla\zeta^T & -\zeta\end{array}\right)
$$
and
$$       Q_2=\left(
\begin{array}{cc} 0 & 0 \\ 0& \frac{\zeta^2+\mu_2(z+1)^2\vert\nabla\zeta\vert^2}{1+\eps_2\zeta}
\end{array}\right).
$$
It follows that
\begin{eqnarray*}
    \nabla^{\mu_2}_{X,z}\cdot Q^{\mu_2}[\eps_2\zeta]\nabla^{\mu_2}_{X,z}
    \underline{\Phi}_{app}
    &=&\dm\Phi^{(0)}\\
    &+ & \eps_2(\dm \cdot \Phi^{(1)} +\nm \cdot
    Q_1\nm\Phi^{(0)})+O(\eps_2^2).
\end{eqnarray*}
Similarly, we obtain
$$
    \partial_n\underline{\Phi}_{app}\,_{\vert_{z=0/-1}}
    =\dz \Phi^{(0)}\,_{\vert_{z=0/-1}}
    +\eps_2\big({\bf e_z}\cdot Q_1\nm\Phi^{(0)}+
    \dz\Phi^{(1)}\big)\,_{\vert_{z=0/-1}}
    +O(\eps_2^2).
$$
Since it is known from Proposition \ref{prop2} that
$$
    \frac{1}{\delta}G^{\mu}[\eps\zeta]\psi_1
    =\sqrt{\mu_2}\nabla\cdot({\mathcal T}_{0,\mu}\nabla\psi_1)
    +\eps_2\mu_2\nabla\cdot (-\zeta+{\mathcal T}_{1,\mu}[\zeta])
    \nabla\psi_1
    +O(\frac{1}{\delta^2}\eps_2^2\mu_2),
$$
one therefore deduces that $\underline{\Phi}_{app}$ solves
(\ref{newbvp}) up to order
$O(\eps_2^2+\frac{1}{\delta^2}\eps_2^2\mu_2)$ provided that
$\Phi^{(0)}$ and $\Phi^{(1)}$ solve
$$
    \left\lbrace
    \begin{array}{l}
    \dm \Phi^{(0)}=0,\\
    \dz \Phi^{(0)}\,_{\vert_{z=0}}=\sqrt{\mu_2}\nabla\cdot ({\mathcal T}_{0,\mu}\nabla\psi_1),\qquad \dz\Phi^{(0)}\,_{\vert_{z=-1}}=0,
    \end{array}\right.
$$
which is obviously solved by $\Phi^{(0)}(X,z)=-\frac{\cosh(\sqrt{\mu_2}(z+1)\vert D\vert)}{\cosh(\sqrt{\mu_2}\vert D\vert)}\frac{\tanh(\sqrt{\mu}\vert D\vert)}{\tanh(\sqrt{\mu_2}\vert D\vert)}\psi_1$, and
$$
    \left\lbrace
    \begin{array}{l}
    \dm \Phi^{(1)}=-\nm\cdot Q_1\nm\Phi^{(0)},\\
    \dz \Phi^{(1)}\,_{\vert_{z=0}}=A,
    \qquad \dz\Phi^{(1)}\,_{\vert_{z=-1}}=0,
    \end{array}\right.
$$
with $A=\mu_2\nabla\cdot (-\zeta+{\mathcal
T}_{1,\mu}[\zeta])\nabla\psi_1-{\bf e_z}\cdot
Q_1\nm\Phi^{(0)}\,_{\vert_{z=0}}$. Because
 $
    -\nm\cdot Q_1\nm \Phi^{(0)}=
    \dm\big[(z+1)\zeta\dz\Phi^{(0)}\big]
$ and
$$
    A=
    \mu_2\nabla\cdot \big[-\zeta+{\mathcal T}_{1,\mu}[\zeta]\nabla\psi_1\big]+\mu_2\nabla\cdot(\zeta\nabla\Phi^{(0)})+
    \dz\big((z+1)\zeta\dz\Phi^{(0)}\big)\,_{\vert_{z=0}},
$$
it results that $\Phi^{(1)}=(z+1)\zeta\dz\Phi^{(0)}+u$, where $u$
solves the boundary-value problem
$$
    \left\lbrace
    \begin{array}{l}
    \dm u=0,\\
    \dz u\,_{\vert_{z=0}}=
    \mu_2\nabla\cdot \big[-\zeta+{\mathcal T}_{1,\mu}[\zeta]\nabla\psi_1
    \big]+\mu_2\nabla\cdot(\zeta\nabla\Phi^{(0)})
    \qquad \dz u\,_{\vert_{z=-1}}=0.
    \end{array}\right.
$$
This latter boundary-value problem can be explicitly solved by
taking the Fourier transform in the horizontal variables and solving
the resulting ordinary differential equation in the variable $z$.
One obtains from this calculation that
\begin{eqnarray*}
    \nabla u_{\vert_{z=0}}&=&\sqrt{\mu_2}\frac{\vert D\vert}{\tanh(\sqrt{\mu_2}\vert D\vert)}\Pi\big[\zeta(1+\frac{\tanh(\sqrt{\mu}\vert D\vert)}{\tanh(\sqrt{\mu_2}\vert D\vert)})\nabla\psi_1\big]\\
&+&
    \sqrt{\mu_2}\nabla\big[
    \frac{\tanh(\sqrt{\mu}\vert D\vert)}{\tanh(\sqrt{\mu_2}\vert D\vert)}
    (\zeta\frac{\tanh(\sqrt{\mu}\vert D\vert)}{\vert D\vert}\Delta\psi_1)\big].
\end{eqnarray*}
Since $\Phi^{(1)}=u+(z+1)\zeta\dz\Phi^{(0)}$ and
$$\nabla [(z+1)\zeta\dz\Phi^{(0)}]_{\vert_{z=0}}=\sqrt{\mu_2}\nabla\big(\zeta\frac{\tanh(\sqrt{\mu}\vert D\vert)}{\vert D\vert}\Delta\psi_1\big),$$
it is deduced immediately that
$
    \nabla\Phi^{(1)}\,_{\vert_{z=0}}=B(\zeta,\nabla\psi_1),
$
where
\begin{eqnarray}\nonumber
    B(\zeta,\nabla\psi_1)&=&
    \sqrt{\mu_2}\frac{\vert D\vert}{\tanh(\sqrt{\mu_2}\vert D\vert)}\Pi\big[\zeta(1+\frac{\tanh(\sqrt{\mu}\vert D\vert)}{\tanh(\sqrt{\mu_2}\vert D\vert)})\nabla\psi_1\big]\\
    \label{exprB}
    &+&
    \sqrt{\mu_2}\nabla\big[
    \big(1+\frac{\tanh(\sqrt{\mu}\vert D\vert)}{\tanh(\sqrt{\mu_2}\vert D\vert)}\big)
    (\zeta\frac{\tanh(\sqrt{\mu}\vert D\vert)}{\vert D\vert}\Delta\psi_1)\big].
\end{eqnarray}
The rigorous result concerning the asymptotic expansion of the
operator  ${\bf H}^{\mu,\delta}[\eps\zeta]$ in the present regime,
which is a corollary of Proposition \ref{expgal}, may now be stated
and proved.
\begin{corollary}[Full dispersion/Full dispersion regime]\label{coro2}
    Let $t_0>d/2$, $s\geq t_0+1/2$, and $\zeta\in H^{s+3/2}(\R^d)$
    be such that
    (\ref{flotte1}) and (\ref{flotte2}) are satisfied. Then, for all
    $\psi_1$ such that $\nabla\psi_1\in H^{s+5/2}(\R^d)$,
    \begin{eqnarray*}
    \big\vert {\bf H}^{\mu,\delta}[\eps\zeta]\psi_1
    -\big(-\frac{\tanh(\sqrt{\mu}\vert D\vert)}{\tanh(\sqrt{\mu_2}\vert D\vert)}\nabla\psi_1+\eps_2 B(\zeta,\nabla\psi_1)\big)\big\vert_{H^s}\\
    \leq \frac{\eps_2^2+\eps^2}{\sqrt{\mu_2}}
    C(\frac{1}{H_1},\frac{1}{H_2},\delta^{max},\mu^{max},
    \mu_2^{max},\vert \zeta\vert_{H^{s+3/2}})
    \vert \nabla\psi_1\vert_{H^{s+5/2}},
    \end{eqnarray*}
    where the bilinear mapping $B(\cdot,\cdot)$ is defined in (\ref{exprB}).
    This estimate is uniform with respect to $\eps\in [0,1]$, $\mu\in (0,\mu^{max})$
and $\delta\in (0,\delta^{max})$ such that
    $\mu_2=\frac{\mu}{\delta^2}\in(0,\mu_2^{max})$.
\end{corollary}
\begin{proof}
The computations above show that
\begin{eqnarray*}
    \nabla^{\mu_2}_{X,z}\cdot Q^{\mu_2}[\eps_2\zeta]\nabla^{\mu_2}_{X,z}
    \underline{\Phi}_{app}
    &=& \eps_2^{2} \nm \cdot {\bf h},
\end{eqnarray*}
with ${\bf h}=Q_1\nm\Phi^{(1)}+Q_2\nm(\Phi^{(0)}+\eps_2\Phi^{(1)})$.
 It is also
easy to check that
\begin{eqnarray*}
    \partial_n\underline{\Phi}_{app}\,_{\vert_{z=0}}&=&
    \sqrt{\mu_2}\nabla\cdot({\mathcal T}_{0,\mu}\nabla\psi_1)
    +\eps_2\mu_2\nabla\cdot (-\zeta+{\mathcal T}_{1,\mu}[\zeta])
    \nabla\psi_1
    + \eps_2^2{\bf e_z}\cdot {\bf h}_{\vert_{z=0}},\\
        \partial_n\underline{\Phi}_{app}\,_{\vert_{z=-1}}&=&
    \eps_2^2{\bf e_z}\cdot {\bf h}_{\vert_{z=-1}}.
\end{eqnarray*}
Therefore, the difference
$v=\underline{\Phi}_{app}-\underline{\Phi}_2$ satisfies  the
boundary-value problem
$$
    \left\lbrace
    \begin{array}{l}
    \nabla^{\mu_2}_{X,z}\cdot Q^{\mu_2}[\eps_2\zeta]\nabla^{\mu_2}_{X,z}v
    =\eps_2^{2} \nm \cdot {\bf h},\\
    \partial_n v\,_{\vert_{z=0}}=
    \sqrt{\mu_2}\nabla\cdot V
    + \eps_2^{2}{\bf e_z}\cdot {\bf h}_{\vert_{z=0}},\qquad
        \partial_nv\,_{\vert_{z=-1}}=
    \eps_2^{2}{\bf e_z}\cdot {\bf h}_{\vert_{z=-1}},
    \end{array}\right.
$$
with $V=    ({\mathcal T}_{0,\mu}\nabla\psi_1)
    +\eps_2\sqrt{\mu_2} (-\zeta+{\mathcal T}_{1,\mu}[\zeta])
    \nabla\psi_1-V^{\mu}[\eps\zeta]\psi_1$, and where
$V^{\mu}[\eps\zeta]\psi_1$ is given by Lemma \ref{lmvert}. Applying
Proposition \ref{expgal} in this situation, it is immediately
deduced that $\vert \nabla v_{\vert_{z=0}}\vert_{H^s}$ is bounded
from above by
$$
    C(\frac{1}{H_2},\delta^{max},\mu_2^{max},\vert \zeta\vert_{H^{s+3/2}})
    \big(\frac{\eps_2^2}{\sqrt{\mu_2}} \Vert {\bf h}\Vert_{H^{s+1/2,1}}
    +\frac{1}{\sqrt{\mu_2}}\vert V \vert_{H^{s+1}}\big).
$$
The stated result is thus a direct consequence of Proposition
\ref{prop2} and the observation that $\Vert {\bf
h}\Vert_{H^{s+1/2,1}}\leq
C(\frac{1}{H_2},\delta^{max},\mu_2^{max},\vert
\zeta\vert_{H^{s+3/2}})\vert\nabla\psi_1\vert_{H^{s+3/2}}$.
\end{proof}

This section concludes with two specializations of Corollary
\ref{coro2} that obtain when additional smallness assumptions are
made on the parameters $\mu$, $\mu_2$ or on $\delta$.  These simple
consequences of Corollary \ref{coro2} will be useful presently.  The
two additional regimes we have in mind are the following.
\begin{enumerate}
    \item The Boussinesq/Full dispersion regime. This regime is
    obtained by assuming that $\mu\sim \eps$ and $\mu_2\sim 1$ (and
    thus $\delta\sim \eps^{1/2}$) in addition to the assumptions
    $\eps\ll 1$ and $\eps_2\ll 1$
    which are required if one wants Corollary \ref{coro2} to provide
    a good approximation.
    \item The Boussinesq/Boussinesq regime. Here, it is assumed
    in addition to $\eps\ll 1$ and $\eps_2\ll1$ that $\mu\sim \eps$ and
    $\mu_2\sim \eps_2$ (and thus $\delta\sim 1$).
\end{enumerate}
\begin{corollary}[Boussinesq/Full dispersion regime]\label{coro2bis}
    Let $t_0>d/2$, $s\geq t_0+1/2$, and $\zeta\in H^{s+3/2}(\R^d)$
    be such that
    (\ref{flotte1}) and (\ref{flotte2}) are satisfied. Then, for all
    $\psi_1$ such that $\nabla\psi_1\in H^{s+5/2}(\R^d)$, the
    inequality
    \begin{eqnarray*}
    \lefteqn{\big\vert {\bf H}^{\mu,\delta}[\eps\zeta]\psi_1
    -\sqrt{\mu}\vert D\vert\coth(\sqrt{\mu_2}\vert D\vert)\big[-\nabla\psi_1-\frac{\mu}{3}\Delta\nabla\psi_1+\eps\Pi\big(\zeta\nabla\psi_1\big)\big]\big\vert_{H^s}}\\
    &\leq& (\frac{\eps_2^2+\eps^2}{\sqrt{\mu_2}}+\eps\mu+\eps\mu^{1/2}\delta)
    C(\frac{1}{H_1},\frac{1}{H_2},\delta^{max},\mu^{max},
    \mu_2^{max},\vert \zeta\vert_{H^{s+3/2}})
    \vert \nabla\psi_1\vert_{H^{s+5/2}},
    \end{eqnarray*}
    where $\Pi=-\frac{\nabla\nabla^T}{\vert D\vert^2}$, holds
     uniformly with respect to $\eps\in [0,1]$, $\mu\in (0,\mu^{max})$
and $\delta\in (0,\delta^{max})$
    such that $\mu_2=\frac{\mu}{\delta^2}\in(0,\mu_2^{max})$.
\end{corollary}
\begin{remark}
    When $\eps\ll 1$, $\mu\sim\eps$, $\mu_2\sim 1$ (and thus $\delta\sim \eps^{1/2}$),
    the three components of the error estimate are all of the same size $O(\eps^2)$.
\end{remark}
\begin{proof}
The result is obtained by using $\tanh(\sqrt{\mu}\vert D\vert)\sim \sqrt{\mu}\vert D\vert-\mu\sqrt{\mu}\frac{1}{3}\D^3$ when $\mu$ is small in Corollary \ref{coro2}.
\end{proof}
Similarly, one may also deduce from Corollary \ref{coro2} the
following result in the Boussinesq-Boussinesq regime.
\begin{corollary}[Boussinesq/Boussinesq regime]\label{coro2ter}
    Let $t_0>d/2$, $s\geq t_0+1/2$, and $\zeta\in H^{s+3/2}(\R^d)$
    be such that
    (\ref{flotte1}) and (\ref{flotte2}) are satisfied. Then, for all
    $\psi_1$ such that $\nabla\psi_1\in H^{s+5/2}(\R^d)$, we have
    \begin{eqnarray*}
    \big\vert {\bf H}^{\mu,\delta}[\eps\zeta]\psi_1
    -\big(-\delta\nabla\psi_1-\frac{\delta}{3}\mu(1-\frac{1}{\delta^2})\Delta\nabla\psi_1+\eps_2(1+\delta)
    \Pi(\zeta\nabla\psi_1)\big)\big\vert_{H^s}\\
    \leq (\frac{\eps_2^2+\eps^2}{\sqrt{\mu_2}}+\mu^{2}+\eps^2)
    C(\frac{1}{H_1},\frac{1}{H_2},\frac{1}{\delta^{min}},\delta^{max},
    \mu^{max},\vert \zeta\vert_{H^{s+3/2}})
    \vert \nabla\psi_1\vert_{H^{s+5/2}},
    \end{eqnarray*}
    where $\Pi=-\frac{\nabla\nabla^T}{\vert D\vert^2}$.
    Moreover, this estimate is uniform with respect to $\eps\in [0,1]$, $\mu\in (0,\mu^{max})$
and $\delta\in (\delta^{min},\delta^{max})$.
\end{corollary}
\begin{remark}
    When $\eps\sim \eps_2 \sim\mu\sim\mu_2 \ll 1$
    (and thus $\delta\sim 1$), the last two
    components of the error estimate are of size $O(\eps^2)$,
    but the first is of size $O(\eps^{3/2})$. This loss of precision
    is not seen at the formal level.  It comes from
    the $1/\sqrt{\mu_2}$ term
    in the elliptic estimate provided by Proposition \ref{expgal}.
\end{remark}
\subsubsection{The Shallow-Water/Shallow-Water Regime: $\mu\ll 1$, $\mu_2\ll1$}
\label{largeeps}

In this regime, large amplitude waves are allowed for the upper
fluid ($\eps=O(1)$) and for the lower fluid ($\eps_2=O(1)$).
Assuming that $\mu\ll1$ and $\mu_2\ll1$ raises the prospect of
making asymptotic expansions of shallow-water type, in terms of
$\mu$ and $\mu_2$. As before, the plan is to formally construct an
approximate solution $\underline{\Phi}_{app}$ to (\ref{newbvp})
having the form
$$
    \underline{\Phi}_{app}={\Phi^{(0)}}+\mu_2{\Phi^{(1)}}.
$$
The formal approximation is then rigorously justified (Corollary
\ref{coro1} below) and the desired expansion results.  From the
expression for $Q^{\mu_2}[\eps_2\zeta]$, we may write
$$
    \nabla^{\mu_2}_{X,z}\cdot Q^{\mu_2}[\eps_2\zeta]\nabla^{\mu_2}_{X,z}
    =\frac{1}{h_2}\dz^2+\mu_2\nabla_{X,z}\cdot Q_1\nabla_{X,z},
$$
with $h_2=1+\eps_2\zeta$ and
$$
Q_1=
    \left(\begin{array}{cc}
    h_2I_{d\times d} & -\eps_2 (z+1)\nabla\zeta\\
     -\eps_2 (z+1)\nabla\zeta^T &
    \frac{\eps_2^2(z+1)^2\vert\nabla\zeta\vert^2}{h_2}
    \end{array}\right).
    $$
It follows readily that
$$
    \nabla^{\mu_2}_{X,z}\cdot Q^{\mu_2}[\eps_2\zeta]\nabla^{\mu_2}_{X,z}
    \underline{\Phi}_{app}
    =\frac{1}{h_2}\dz^2\Phi^{(0)}
    +\mu_2\big(\nabla_{X,z}\cdot Q_1\nabla_{X,z}\Phi^{(0)}+
    \frac{1}{h_2}\dz^2\Phi^{(1)}\big)
    +O(\mu_2^2).
$$
Similarly, one infers that at $z=0$ and $z=-1$,
$$
    \partial_n\underline{\Phi}_{app}
    =\frac{1}{h_2}\dz \Phi^{(0)}
    +\mu_2\big({\bf e_z}\cdot Q_1\nabla_{X,z}\Phi^{(0)}+
    \frac{1}{h_2}\dz\Phi^{(1)}\big)
    +O(\mu_2^2).
$$
Since it is known from Proposition \ref{prop1} that
$$
    \frac{1}{\delta}G^{\mu}[\eps\zeta]\psi_1=
    \delta\mu_2\nabla\cdot(h_1\nabla\psi_1)+O(\frac{\mu^2}{\delta})
$$
(with $h_1=1-\eps\zeta$), it is clearly the case that
$\underline{\Phi}_{app}$ solves (\ref{newbvp}) up to order
$O(\mu_2^2+\frac{\mu^2}{\delta})$ provided that $\Phi^{(0)}$ and
$\Phi^{(1)}$ solve
$$
    \left\lbrace
    \begin{array}{l}
    \dz^2 \Phi^{(0)}=0,\\
    \dz \Phi^{(0)}\,_{\vert_{z=0}}=0,\qquad \dz\Phi^{(0)}\,_{\vert_{z=-1}}=0,
    \end{array}\right.
$$
which is obviously solved by any $\Phi^{(0)}(X,z)=\Phi^{(0)}(X)$ independent of $z$, and
$$
    \left\lbrace
    \begin{array}{l}
    \dz^2 \Phi^{(1)}=-h_2^2\Delta\Phi^{(0)},\\
    \dz \Phi^{(1)}\,_{\vert_{z=0}}=h_2\big(\eps_2\nabla\zeta\cdot\nabla\Phi^{(0)}+\delta\nabla\cdot(h_1\nabla\psi_1)\big),
    \qquad \dz\Phi^{(1)}\,_{\vert_{z=-1}}=0,
    \end{array}\right.
$$
where we have used the fact that $\Phi^{(0)}$ does not depend on
$z$. Solving this second order ordinary differential equation in the
variable $z$ with the boundary condition at $z=0$ yields (up to a
function independent of $z$ which we take equal to $0$ for the sake
of simplicity),
$$
    \Phi^{(1)}=-\frac{z^2}{2}h_2^2\Delta\Phi^{(0)}+z(\dz \Phi_1\,_{\vert_{z=0}}).
$$
Matching the boundary condition at $z=-1$ leads to the restriction
$$
    \nabla\cdot(h_2\nabla\Phi^{(0)})=-\delta\nabla\cdot(h_1\nabla\psi_1),
$$
    which implies that
$\Pi (h_2\nabla\Phi^{(0)})=\Pi(-\delta h_1\nabla\psi_1),$ where
$\Pi=-\frac{\nabla\nabla^T}{\vert D\vert^2}$ is the orthogonal
projector onto the gradient vector fields of $L^2(\R^d)^d$ defined
earlier. We will solve this equation thanks to the following lemma.
\begin{lemma}   \label{lemmproj}
Assume that $\zeta\in L^\infty(\R^d)$ is such that $\vert\eps_2\zeta\vert_\infty<1$. 
 Let also
$W\in L^2(\R^d)^d$. Then\\
i. One can define the mapping ${\mathfrak Q}[\eps_2\zeta]$ as
$$
	{\mathfrak Q}[\eps_2\zeta]:\qquad
	\begin{array}{ccc}
	\displaystyle L^2(\R^d)^d& \to &L^2(\R^d)^d\\
	\displaystyle U&\mapsto & \displaystyle \sum_{n=0}^\infty (-1)^n(\Pi(\eps_2\zeta\Pi\cdot))^n (\Pi U)
	\end{array}
$$
ii. There exists a unique solution $V\in L^2(\R^d)^d$ to the equation
$$
	\nabla\cdot(h_2 V)=\nabla\cdot W,\qquad (h_2=1+\eps_2\zeta)
$$
such that $\Pi V=V$ and one has $V={\mathfrak Q}[\eps_2\zeta] W$;\\
iii- If moreover $\zeta\in H^s(\R^d)$ and $W\in H^s(\R^d)^d$ ($s>d/2+1$) then
${\mathfrak Q}[\eps_2\zeta]W\in H^s(\R^d)^d$ and
$$
	\vert {\mathfrak Q}[\eps_2\zeta]W\vert_{H^s}
	\leq C(\vert \eps_2\zeta\vert_{H^s},\frac{1}{1-\vert\eps_2\zeta\vert_\infty})
	\vert W\vert_{H^s}.
$$
\end{lemma}
\begin{remark}\label{remafix1}
In dimension $d=1$, one has $\Pi=1$ and the first point of the lemma simplifies
into $V=\frac{1}{h_2}W$ so that the proof is trivial.
\end{remark}
\begin{proof}
  i. The result follows from the observation that under the assumptions of the lemma,
one has 
\begin{equation}\label{boston1}
	\Vert \Pi(\eps_2\zeta\Pi\cdot)\Vert_{L^2\to L^2}\leq \vert\eps_2\zeta\vert_\infty
	<1,
\end{equation}
so that the series used to define ${\mathfrak Q}[\eps_2\zeta]U$ converges in $L^2(\R^d)^d$.\\
ii. Let us first check that $V={\mathfrak Q}[\eps_2\zeta]W$ is indeed a solution of the equation
stated in the lemma.  Since $V=\Pi V$, one can remark that 
\begin{eqnarray*}
	\nabla\cdot (\eps_2\zeta V)&=&\nabla\cdot ( \Pi(\eps_2\zeta\Pi V))\\
		&=&- \nabla\cdot \sum_{n=1}^\infty (-1)^n(\Pi(\eps_2\zeta\Pi\cdot))^n(\Pi W)\\
		&=&- \nabla\cdot (V-\Pi W),
\end{eqnarray*}
from which one deduces easily that $\nabla\cdot (h_2 V)=\nabla\cdot W$.\\
Let us now turn to prove uniqueness of the solution by proving that one has necessarily 
$V=0$ if $W$=0. To check that this is the case, just remark that from the equation
$\nabla\cdot (h_2V)=0$ and the requirement that $\Pi V=V$, one has
$$
	V=-\Pi(\eps_2\zeta\Pi V);
$$
since $\Vert \Pi(\eps_2\zeta \Pi V)\Vert_{L^2\to L^2}\leq \vert\eps_2\zeta\vert_\infty<1$,
it follows easily  that $V=0$.\\
iii. It is clear from (\ref{boston1}) that $\vert {\mathfrak Q}[\eps_2\zeta]W\vert_2\leq\frac{1}{1-\vert\eps_2\zeta\vert_\infty}$. Now, applying
$\Lambda^s$ to the equations, one gets
$$
	\nabla\cdot (h_2\Lambda^s V)=\nabla\cdot \widetilde{W},
$$
with $\widetilde{W}=\Lambda^s W+[\Lambda^s ,\eps_2\zeta] V$. The result follows therefore
from the $L^2$ estimate, a standard commutator estimate and a simple induction.  
\end{proof}
If Lemma \ref{lemmproj} is applied with $V=\nabla\Phi^{(0)}$, $W=-\delta
h_1\nabla\psi_1$, there results the
equation
$$
    \nabla \Phi^{(0)}=-\delta\ {\mathfrak Q}[\eps_2\zeta](h_1\nabla\psi_1).
$$
Note that when d = 1, this reduces to
$$
  \dx \Phi^{(0)} = -\delta\frac{h_1}{h_2}\dx\psi_1.
$$
The following corollary of Proposition \ref{expgal}, which gives the
needed asymptotic expansion of the operator ${\bf
H}^{\mu,\delta}[\eps\zeta]$ in the present regime, now comes into
view.
\begin{corollary}[Shallow water/Shallow water regime]\label{coro1}
    Let $t_0>d/2$, $s\geq t_0+1/2$, and $\zeta\in H^{s+3/2}(\R^d)$
    be such that
    (\ref{flotte1}) and (\ref{flotte2}) are satisfied. Let $h_1=1-\eps\zeta$
    and $h_2=1+\eps_2\zeta$ and let
    $\psi_1$ be such that $\nabla\psi_1\in H^{s+5/2}(\R^d)$.  Then
    it follows that
    \begin{eqnarray*}
    \lefteqn{\vert {\bf H}^{\mu,\delta}[\eps\zeta]\psi_1+\delta\ {\mathfrak Q}[\eps_2\zeta](h_1\nabla\psi_1)\vert_{H^s}}\\
    &\leq& \delta(\mu+\mu_2)
    C\big((1-\delta(1-H_1))^{-1},\frac{1}{H_2},\eps_2^{max},
    \mu_2^{max},\vert \zeta\vert_{H^{s+3/2}}\big)
    \vert \nabla\psi_1\vert_{H^{s+5/2}},
    \end{eqnarray*}
uniformly with respect to $\eps\in [0,1]$, $\mu\in (0,1)$
and $\delta<\frac{1}{1-H_1}$ such that $\eps_2=\eps\delta\in[0,\eps_2^{max}]$ and $\mu_2=\frac{\mu}{\delta^2}\in(0,\mu_2^{max})$.
\end{corollary}
\begin{remark}
    When $\eps\sim\eps_2\sim\mu\sim \mu_2 \ll 1$ (and thus $\delta\sim 1$),
    one deduces from the above corollary that
    $ {\bf H}^{\mu,\delta}[\eps\zeta]\psi_1=-\delta \nabla\psi_1
    +O(\eps)$, which
    is consistent with the asymptotic expansion provided by
    Corollary \ref{coro2ter}. A similar matching would have been
    observed for the next order terms if we had computed them
     in Corollary \ref{coro1}.
\end{remark}
\begin{remark}
When d = 1, one has $\delta\ {\mathfrak Q}[\eps_2\zeta](h_1\nabla\psi_1)=
\delta\frac{h_1}{h_2}\dx\psi_1$.
\end{remark}
\begin{proof}
Since (\ref{flotte1}), (\ref{flotte2}) and the condition $\delta(1-H_1)<1$ imply that
$\vert\eps_2\zeta\vert_\infty<1$, one can use Lemma \ref{lemmproj}  and 
the computations above indicate that
\begin{eqnarray*}
    \nabla^{\mu_2}_{X,z}\cdot Q^{\mu_2}[\eps_2\zeta]\nabla^{\mu_2}_{X,z}
    \underline{\Phi}_{app}&=&
    \mu_2^2 \nabla_{X,z} \cdot Q_1\nabla_{X,z}\Phi_1\\
    &=& \mu_2^{3/2} \nm \cdot {\bf h},
\end{eqnarray*}
with
$${\bf h}=\left(\begin{array}{cc} I_{d\times d}& 0\\ 0 & \sqrt{\mu_2}
            \end{array}\right)Q_1\nabla_{X,z}\Phi_1.
            $$
            It is also easy to check that
\begin{eqnarray*}
    \partial_n\underline{\Phi}_{app}\,_{\vert_{z=0}}&=&
    \delta\mu_2\nabla\cdot(h_1\nabla\psi_1)
    + \mu_2^{3/2}{\bf e_z}\cdot {\bf h}_{\vert_{z=0}},\\
        \partial_n\underline{\Phi}_{app}\,_{\vert_{z=-1}}&=&
    \mu_2^{3/2}{\bf e_z}\cdot {\bf h}_{\vert_{z=-1}}.
\end{eqnarray*}
Thus, the difference $u=\underline{\Phi}_{app}-\underline{\Phi}_2$
satisfies  the boundary-value problem
$$
    \left\lbrace
    \begin{array}{l}
    \nabla^{\mu_2}_{X,z}\cdot Q^{\mu_2}[\eps_2\zeta]\nabla^{\mu_2}_{X,z}u
    =\mu_2^{3/2} \nm \cdot {\bf h},\\
    \partial_n u\,_{\vert_{z=0}}=
    \sqrt{\mu_2}\nabla\cdot V
    + \mu_2^{3/2}{\bf e_z}\cdot {\bf h}_{\vert_{z=0}},\qquad
        \partial_nu\,_{\vert_{z=-1}}=
    \mu_2^{3/2}{\bf e_z}\cdot {\bf h}_{\vert_{z=-1}}.
    \end{array}\right.
$$
with $V=h_1\delta\sqrt{\mu_2}\nabla\psi_1-V^{\mu}[\eps\zeta]\psi_1$,
where $V^{\mu}[\eps\zeta]\psi_1$ is given by Lemma \ref{lmvert}.
One concludes from Proposition \ref{expgal}  that $\vert \nabla
u_{\vert_{z=0}}\vert_{H^s}$ is bounded from above by the quantity
$$
    C\Big(\frac{1}{H_2},\eps_2^{max},\mu_2^{max},\vert \zeta\vert_{H^{s+3/2}}\Big)
    \Big(\mu_2 \Vert {\bf h}\Vert_{H^{s+1/2,1}}
    +\delta\vert
    h_1\nabla\psi_1-\frac{1}{\sqrt{\mu}}V^{\mu}[\eps\zeta]\psi_1\vert_{H^{s+1}}\Big).
$$
The result is a direct consequence of Proposition \ref{prop1} since
$$\Vert {\bf h}\Vert_{H^{s+1/2,1}}\leq \delta
C(\frac{1}{H_2},\eps_2^{max},\mu_2^{max},\vert
\zeta\vert_{H^{s+3/2}})\vert\nabla\psi_1\vert_{H^{s+3/2}}$$.
\end{proof}

\subsubsection{The Shallow-Water/Small-Amplitude Regime: $\mu\!\ll\! 1$, $\eps_2\!\ll\!1$}\label{sectSWSA}
It is now presumed that both $\mu$ and $\eps_2$ are small, but no
such restriction is laid upon $\eps$ nor $\mu_2$.  So, this regime
is not a subcase of the regimes investigated in Sections
\ref{smalleps} and \ref{largeeps}. We construct an approximate
solution $\underline{\Phi}_{app}$ to (\ref{newbvp}) exactly as in \S
\ref{smalleps}, but only a first-order approximation of the form
$$
    \underline{\Phi}_{app}={\Phi^{(0)}};
$$
will be required.  Since $\mu\ll1$ here,  Proposition \ref{prop2}
may be utilized to write
$$
    \frac{1}{\delta}G^\mu[\eps\zeta]\psi_1=\frac{\mu}{\delta}
    \nabla\cdot\big(h_1\nabla\psi_1\big)+O(\frac{\mu^2}{\delta}).
$$
Just as in \S \ref{smalleps}, it can be shown that $\Phi^{(0)}$ must
solve the boundary-value problem
$$
    \left\lbrace
    \begin{array}{l}
    \dm \Phi^{(0)}=0,\\
    \dz \Phi^{(0)}\,_{\vert_{z=0}}=
    \frac{\mu}{\delta}
    \nabla\cdot\big(h_1\nabla\psi_1\big)
    ,\qquad \dz\Phi^{(0)}\,_{\vert_{z=-1}}=0,
    \end{array}\right.
$$
which is to say that
$$
    \Phi^{(0)}(X,z)=\sqrt{\mu}\frac{\cosh(\sqrt{\mu_2}(z+1)\vert D\vert)}{\cosh(\sqrt{\mu_2}\vert D\vert)}\frac{1}{\D\tanh(\sqrt{\mu_2}\vert D\vert)}
    \nabla\cdot\big(h_1\nabla \psi_1\big).
$$
The following result is proved exactly as was Corollary \ref{coro2}.
\begin{corollary}[Shallow water/Small amplitude regime]\label{coro3}
    Let $t_0>d/2$, $s\geq t_0+1/2$, and $\zeta\in H^{s+3/2}(\R^d)$
    be such that
    (\ref{flotte1}) and (\ref{flotte2}) are satisfied. Then, for all
    $\psi_1$ such that $\nabla\psi_1\in H^{s+5/2}(\R^d)$,
    it is the case that
    \begin{eqnarray*}
    \lefteqn{\big\vert {\bf H}^{\mu,\delta}[\eps\zeta]\psi_1
    +\sqrt{\mu}\vert D\vert\coth(\sqrt{\mu_2}\D)\Pi
    \big(h_1\nabla\psi_1\big)\big\vert_{H^s}}\\
    &\leq& \frac{\mu^{3/2}+\eps_2\sqrt{\mu}}{\sqrt{\mu_2}}
    C(\frac{1}{H_1},\frac{1}{H_2},\delta^{max},
    \mu_2^{max},\vert \zeta\vert_{H^{s+3/2}})
    \vert \nabla\psi_1\vert_{H^{s+5/2}},
    \end{eqnarray*}
    where $h_1=1-\eps\zeta$ and $\Pi=-\frac{\nabla\nabla^T}{\D^2}$
    is given by (\ref{exprB}). This estimate is
    uniform with respect to $\eps\in [0,1]$, $\mu\in (0,1)$
    and $\delta\in (0,\delta^{max})$ such that
    $\mu_2=\frac{\mu}{\delta^2}\in(0,\mu_2^{max})$.
\end{corollary}
\begin{remark}\label{remaILWBO}
    Several regimes fall within the range of Corollary \ref{coro3}.
    \begin{itemize}
    \item\label{remit1}
    The SW/FD regime: when $\mu\ll 1$, $\eps_2\ll 1$ and $\eps\sim \mu_2\sim 1$
    (and thus $\delta^2\sim \mu\sim \eps_2^2$);
    the precision of the approximation is $O(\mu)$.
    \item\label{remit2} The ILW regime: if
    $\mu\sim \eps^2\ll 1$ and $\mu_2\sim 1$
    (and thus $\delta^2\sim \mu\sim\eps_2$); in this case, the estimate in the corollary
    can be simplified without adverse effects on the precision of the
    approximation to simply
    \begin{equation}\label{asympILW}
    {\bf H}^{\mu,\delta}[\eps\zeta]\psi_1
    =-\sqrt{\mu}\vert D\vert\coth(\sqrt{\mu_2}\D)
    \nabla\psi_1+O(\mu).
    \end{equation}
    \item\label{remit3}The BO regime: if $\mu\ll 1$ and $\delta=0$ (and thus
    $\mu_2=\infty$, $\eps_2=0$), one gets formally from (\ref{asympILW})
    that
    \begin{equation}\label{asympBO}
    {\bf H}^{\mu,\delta}[\eps\zeta]\psi_1
    \sim-\sqrt{\mu}\vert D\vert
    \nabla\psi_1.
    \end{equation}
    \end{itemize}
\end{remark}

\section{Asymptotic Models for Internal Waves}

The preliminary analysis in Section 2 allows us to derive the
various asymptotic models referred to in the Introduction.

\subsection{The small amplitude/small amplitude regime: $\eps\ll 1$, $\eps_2\ll1$}\label{SA-SA}

Derived first are various models corresponding to the case wherein
the interface deformation are small for both fluids. Different
systems of equations obtain, depending on the sizes of the
paramaters $\eps$, $\mu$ and $\delta$ (and thus $\eps_2$ and
$\mu_2$).

\subsubsection{The Full Dispersion/Full Dispersion Regime: $\eps\sim \eps_2\ll 1$ and $\mu\sim \mu_2=O(1)$}\label{ILW-ILW}

An asymptotic model can be derived from (\ref{eqdepndim}) by
replacing the operators $G^\mu[\eps\zeta]$ and ${\bf
H}^{\mu,\delta}[\eps\zeta]$ by their asymptotic expansions, provided
by Proposition \ref{prop2} and Corollary \ref{coro2} in the present
regime. The following theorem shows that in the present regime, the
internal wave equations are consistent  with following FD/FD system;
\begin{equation}\label{eqFDFD}
    \left\lbrace
    \begin{array}{l}
    \dt \zeta+
    \frac{1}{\sqrt{\mu}}\frac{\nabla}{\D}\cdot \big(\frac{\tm\tmd}{\gamma \tmd+\tm}{\bf v}\big)
    \\
    \indent
    +\frac{\eps_2}{\sqrt{\mu}}\frac{\nabla}{\D}\cdot \big(\frac{\tm\tmd}{\gamma \tmd+\tm}B(\zeta,\frac{\tmd}{\gamma\tmd+\tm}{\bf v})\big)\\
    \indent
    -\eps\nabla\cdot(\zeta\frac{\tmd}{\gamma\tmd+\tm}{\bf v} )
    +\eps\D\tm\big(\zeta\frac{\nabla}{\D}\cdot(\frac{\tm\tmd}{\gamma\tmd+\tm}{\bf v})\big)=0\\
    \dt {\bf v}
    +(1-\gamma)\nabla\zeta\\
    \indent+\frac{\eps}{2}\nabla\big(
    \big\vert \frac{\tm}{\gamma\tmd+\tm}{\bf v}\big\vert^2-\gamma
    \big\vert \frac{\tmd}{\gamma\tmd+\tm}{\bf v}\big\vert^2
    \big)
    +\eps\frac{\gamma-1}{2}\nabla \big(\frac{\nabla}{\D}\cdot(\frac{\tm\tmd}{\gamma\tmd+\tm}{\bf v})\big)^2=0,
    \end{array}\right.
\end{equation}
where as before, $\tm=\tanh(\sqrt{\mu}\D)$,
$\tmd=\tanh(\sqrt{\mu_2}\D)$ and the bilinear mapping
$B(\cdot,\cdot)$ is given by (\ref{exprB}).
\begin{theorem}\label{propFDFD}
    Let $0<\delta^{min}<\delta^{max}$.
    The internal waves equations
    (\ref{eqdepndim}) are consistent with the FD/FD equations
    (\ref{eqFDFD}) in the sense of Definition \ref{defcons},
    with a precision $O(\eps^2)$, and uniformly with respect
    to
    $\eps\in[0,1]$, $\mu\in (0,\mu^{max})$ and
    $\delta\in [\delta^{min},\delta^{max}]$.
\end{theorem}
\begin{remark}
    One can give a more precise
    estimate of the precision, as in Corollary \ref{coro2}
    for instance. It simplifies the exposition to use the
    notation $O(\eps^2)$ and the associated rough estimate
    of the precision.  We follow this policy throughout
    the discussion.
\end{remark}
\begin{remark}
    It is straightforward to check that the dispersion relation of (\ref{eqFDFD})
    is exactly the same as (\ref{disprel}), which
    is the reason we refer to (\ref{eqFDFD}) as
    a ``full dispersion'' model. In particular, (\ref{eqFDFD}) is
    linearly well-posed provided that $\gamma<1$.
\end{remark}

\begin{proof}
First, notice that with the range of parameters considered in the
theorem, one has $\eps\sim\eps_2$ when $\eps\to 0$, while $\mu\sim
\mu_2=O(1)$. By the definition (\ref{defv}) of ${\bf v}$ and using
Proposition \ref{prop2} and Corollary \ref{coro2}, one deduces from
(\ref{eqdepndim}) that
$$
    \left\lbrace
    \begin{array}{l}
    \dt \zeta-\frac{1}{\sqrt{\mu}}\frac{\nabla}{\D}\cdot(\tm
    \nabla\psi_1)
    +\eps\nabla\cdot(\zeta\nabla\psi_1)
    -\eps\vert D\vert \tm
    \big(\zeta\frac{\nabla}{\vert D\vert}\cdot(\tm\nabla\psi_1)\big)
    =O(\eps^2)
    \\
    \dt {\bf v}
    +(1-\gamma)\nabla\zeta\\
    \indent\indent+\frac{\eps}{2}\nabla(
    \vert {\bf H}^{\mu,\delta}[\eps\zeta]\psi_1\vert^2
    -\gamma\vert \nabla\psi_1\vert^2)+\eps\frac{\gamma-1}{2}\nabla\big(\frac{\nabla}{\vert D\vert}\cdot(\tm\nabla\psi_1)\big)^2=O(\eps^2).
    \end{array}\right.
$$
It follows from Corollary \ref{coro2} and the relation
 ${\bf H}^{\mu,\delta}[\eps\zeta]\psi_1={\bf v}+\gamma\nabla\psi_1$
that
$$
    \nabla\psi_1=-
    \frac{\tmd}{\gamma\tmd+\tm}\big({\bf v}
    +\eps_2
    B(\zeta, \frac{\tmd}{\gamma\tmd+\tm}{\bf v})\big)
    +O(\eps^2).
$$
The result is view is now apparent.
\end{proof}

\subsubsection{The Boussinesq/Full Dispersion Regime $\mu\sim\eps\ll 1$,
$\mu_2\sim 1$}\label{B-ILW}

We show here that in this regime (for which one also has
$\delta^2\sim\eps$ and thus $\eps_2\sim \eps^{3/2}\ll 1$), the
internal waves equations (\ref{eqdepndim}) are consistent with the
\emph{three-parameter family} of Boussinesq/FD systems
\begin{equation}\label{eqB-FD}
          \left\lbrace
                  \begin{array}{l}
                  \dsp \big(1-\mu b\Delta\big)\dt \zeta+
                  \frac{1}{\gamma}\nabla\cdot\big((1-\eps\zeta){\bf v}_\beta\big)\vspace{1mm}\\
    \dsp \indent-\frac{\sqrt{\mu}}{\gamma^2}\D\coth(\sqrt{\mu_2}\D)\nabla\cdot{\bf v}_\beta
+\frac{\mu}{\gamma}\Big(a-\frac{1}{\gamma^2}\coth^2(\sqrt{\mu_2}\D)\Big)\Delta\nabla\cdot{\bf
v}_\beta = 0\vspace{1mm}\\
                \dsp(1-\mu d\Delta)\dt\mathbf{v}_\beta +(1-\gamma)\nabla\zeta-\frac{\eps}{2\gamma}\nabla\vert{\bf v}_\beta\vert^2
    +\mu c(1-\gamma)\Delta\nabla\zeta=0,
            \end{array}\right.
\end{equation}
where ${\bf v}_\beta=(1-\mu\beta \Delta)^{-1}{\bf v}$ and the
constants $a$, $b$, $c$ and $d$ are defined now.
\begin{theorem}\label{propBFD}
    Let $0<c^{min}<c^{max}$, $0<\mu^{min}_2<\mu^{max}_2$,
     and set
    $$
    a=\frac{1}{3}(1-\alpha_1-3\beta),\qquad
    b=\frac{1}{3}\alpha_1,\qquad
    c=\beta\alpha_2,\qquad
    d=\beta(1-\alpha_2),
    $$
    with $\alpha_1\geq 0$, $\beta\geq 0$ and $\alpha_2\leq 1$.
    With these choices of parameters, the internal wave equations
    (\ref{eqdepndim}) are consistent with the Boussinesq/FD equations
    (\ref{eqB-FD}) in the sense of Definition \ref{defcons},
    with a precision $O(\eps^{3/2})$, and uniformly with respect
    to
    $\eps\in[0,1]$, $\mu\in (0,1)$ and $\delta\in (0,1)$ satisfying
    the conditions
    $$
    c^{min}\leq \frac{\eps}{\mu}\leq c^{max}
    \quad\mbox{ and }\quad
    \mu_2^{min}\leq \frac{\mu}{\delta^2}\leq \mu_2^{max}.
    $$
\end{theorem}
\begin{remark}
The dispersion relation associated to (\ref{eqB-FD}) is
$$
    \omega^2=\frac{1-\gamma}{\gamma}\vert {\bf k}\vert^2(1-\mu c\vert{\bf k}\vert^2)\frac{1-\frac{\sqrt{\mu}}{\gamma}\vert {\bf k}\vert\coth(\sqrt{\mu_2} \vert {\bf k}\vert)-\mu \vert {\bf k}\vert^2\big(a-\frac{1}{\gamma^2}\coth^2(\sqrt{\mu_2}\vert {\bf k}\vert)\big)}{(1+\mu b\vert {\bf k}\vert^2)(1+\mu d\vert {\bf k}\vert^2)},
$$
and (\ref{eqB-FD}) is therefore linearly well-posed when $b,d\geq 0$
and $a,c\leq 0$. Notice that in the case
$\alpha_1=\alpha_2=\beta=0$, one has $a=\frac{1}{3}$ and $b=c=d=0$
and the corresponding system is thus linearly ill-posed. The freedom
to choose a well-posed model is just one of the advantages of a
three-parameter family of formally equivalent systems. The same
remark has already been made about  the Boussinesq systems for wave
propagation in the case of surface gravity waves \cite{BCS,BCL}).
\end{remark}
\begin{proof}  The proof is made in several steps, corresponding to
particular assumptions about the parameters $\alpha_1, \alpha_2$ and
$\beta$.  Throughout, use will be freely made of the relations
$\mu\sim\eps$ and $\mu_2\sim 1$.\\

 {\bf Step 1.} {\emph The case $\alpha_1=0$, $\beta=0$,
$\alpha_2=0$.}  From the expansion of the Dirichlet-Neumann
operator, see Remark \ref{remB}, it follows as in the previous
section that
$$
\left\lbrace
    \begin{array}{l}
    \dt \zeta-\nabla\cdot ((1-\eps\zeta)\nabla\psi_1\ )-\frac{\mu}{3}\nabla\cdot\Delta\nabla\psi_1=O(\eps^2)
    \\
    \dt {\bf v}
    +(1-\gamma)\nabla\zeta +\frac{\eps}{2}\nabla(
    \vert {\bf H}^{\mu,\delta}[\eps\zeta]\psi_1\vert^2
    -\gamma\vert \nabla\psi_1\vert^2)=O(\eps^2),
    \end{array}\right.
$$
where the fact that $O(\mu)=O(\eps)$ has been used. From the
relation  ${\bf H}^{\mu,\delta}[\eps\zeta]\psi_1={\bf
v}+\gamma\nabla\psi_1$, and Corollary \ref{coro2bis}, it is seen
that
$$
\nabla\psi_1=-\frac{1}{\gamma}{\bf v}-\frac{\sqrt{\mu}}{\gamma}\frac{\D}{\tmd}\big[1+\frac{\mu}{3}\Delta
 -\eps  \Pi(\zeta\cdot)\big]\nabla\psi_1+O(\eps^2).
$$
Again using the fact that $O(\mu)=O(\eps)$, one concludes that
$$
\nabla\psi_1=-\frac{1}{\gamma}{\bf v}+\frac{\sqrt{\mu}}{\gamma^2}\frac{\D}{\tmd}{\bf v}+\frac{\mu}{\gamma^3}\frac{\Delta}{\tmd^2}{\bf v}+O(\eps^\frac{3}{2})
$$
and the result follows.\\
{\bf Step 2.} {\emph The case $\alpha_1 \geq 0$, $\beta=0$,
$\alpha_2=0$.} We use here the the classical BBM trick \cite{BBM}.
It is clear from the first equation that
$$
    \dt \zeta=-\frac{1}{\gamma}\nabla\cdot {\bf v}+O(\eps^{1/2}),
$$
from which it is inferred that
$$
    \nabla\cdot {\bf v}=(1-\alpha_1)\nabla\cdot {\bf v}-\alpha_1\gamma
    \dt\zeta+O(\eps^{1/2}).
$$
Replacing $\nabla\cdot{\bf v}$ by this expression in the component
$\frac{\mu}{3\gamma}\Delta\nabla\cdot{\bf v}$ of the first equation of the
system derived in Step 1, leads to the desired result.\\
{\bf Step 3.} {\emph The case $\alpha_1\geq 0$, $\beta\geq 0$,
$\alpha_2=0$.} Replacing ${\bf v}$ by $(1-\mu\beta\Delta){\bf
v}_\beta$ in the system of equations derived in Step 2, and
neglecting the $O(\eps^{3/2})$ terms is all that is required in this
case.\\
{\bf Step 4.} {\emph The case  $\alpha_1\geq 0$, $\beta\geq 0$,
$\alpha_2\leq 1$.} We use once again the BBM trick.  From the second
equation in the system derived in Step 3, one obtains that for all
$\alpha_2\leq 1$,
$$
    \dt {\bf v}_\beta=(1-\alpha_2)\dt {\bf
    v}_\beta-\alpha_2(1-\gamma)\nabla\zeta+O(\eps).
$$
If this relationship is substituted into the system derived in Step
3, the result follows.
\end{proof}

\subsubsection{The Boussinesq/Boussinesq Regime $\eps\sim\mu\sim \eps_2\sim\mu_2\ll 1$}\label{B-B}

In this regime, the nonlinear and dispersive effects are of the same
size for both fluids; the systems of equations that are derived from
the internal waves equations (\ref{eqdepndim}) in this situation are
the following \emph{three-parameter family} of Boussinesq/Boussinesq
systems, {\it viz.}
\begin{equation}\label{eqB-B}
          \left\lbrace
                  \begin{array}{l}
                 \dsp  \Big(1-\mu b\Delta\Big)\dt \zeta+\frac{1}{\gamma+\delta}\nabla\cdot{\bf v}_\beta
    +\eps\frac{\delta^2-\gamma}{(\gamma+\delta)^2}
    \nabla\cdot(\zeta{\bf v}_\beta)
        +\mu a
    \nabla\cdot\Delta {\bf v}_\beta=0\vspace{1mm}
                 \\
                \dsp \Big(1-\mu d\Delta\Big)\dt\mathbf{v}_\beta +(1-\gamma)\nabla\zeta+\frac{\eps}{2}\frac{\delta^2-\gamma}{(\delta+\gamma)^2}\nabla\vert{\bf v}_\beta\vert^2 +\mu c\Delta\nabla\zeta=0,
            \end{array}\right.
\end{equation}
where ${\bf v}_\beta=(1-\mu\beta\Delta)^{-1}{\bf v}$, and where the
coefficients $a,b,c,d$ are provided in the statement of the next
theorem.
\begin{theorem}\label{propBB}
    Let $0<c^{min}<c^{max}$, $0<\delta^{min}<\delta^{max}$,
     and set
    $$
    \begin{array}{ll}
    \dsp a=\frac{(1-\alpha_1)(1+\gamma\delta)-3\delta\beta(\gamma+\delta)}{3\delta(\gamma+\delta)^2},&
    b=\alpha_1\frac{1+\gamma\delta}{3\delta(\gamma+\delta)},\\
    c=\beta\alpha_2,&
    d=\beta(1-\alpha_2),
    \end{array}
    $$
    with $\alpha_1\geq 0$, $\beta\geq 0$ and $\alpha_2\leq 1$.
    With this specification of the parameters,
    The internal wave equations
    (\ref{eqdepndim}) are consistent with the Boussinesq/Boussinesq
    equations
    (\ref{eqB-B}) in the sense of Definition \ref{defcons},
    with a precision $O(\eps^2)$, and uniformly with respect
    to
    $\eps\in[0,1]$, $\mu\in (0,1)$ and
    $\delta\in [\delta^{min},\delta^{max}]$ such that
    $c^{min}<\frac{\eps}{\mu}<c^{max}$.
\end{theorem}
\begin{remark}
Taking $\gamma=0$ and $\delta=1$ in the Boussinesq/Boussinesq
equations (\ref{eqB-B}),
 reduces them to the system
$$
 \left\lbrace
                 \begin{array}{l}
         \dsp \big(1-\mu\frac{\alpha_1}{3}\Delta\big)\dt \zeta+\nabla\cdot((1+\eps \zeta){\bf v})
        +\mu\frac{1-\alpha_1-3\beta}{3}
    \nabla\cdot\Delta {\bf v} =0\vspace{1mm}\\
    \dsp \big(1-\mu\beta(1-\alpha_2)\Delta\big)
    \dt\mathbf{v} +\nabla\zeta+\frac{\eps}{2}\nabla\vert{\bf v}\vert^2
    +\mu\beta\alpha_2\Delta\nabla\zeta =0,
            \end{array}\right.
$$
which is exactly the family of formally equivalent Boussinesq
systems derived in \cite{BCS,BCL}.
\end{remark}
\begin{remark}
The dispersion relation associated to  (\ref{eqB-B}) is
$$
    \omega^2=\vert{\bf k}\vert^2\frac{(\frac{1}{\gamma+\delta}
    -\mu a \vert {\bf k}\vert^2)(1-\gamma-\mu c\vert{\bf k}\vert^2)}{(1+\mu b\vert {\bf k}\vert^2)(1+\mu d \vert {\bf
    k}\vert^2)}.
$$
It follows that (\ref{eqB-B}) is linearly well-posed when $a,c\leq
0$ and $b,d\geq 0$. The system corresponding to
$\alpha_1=\alpha_2=\beta=0$ is ill-posed (one can check that
$a=\frac{1+\gamma\delta}{3\delta(\gamma+\delta)^2}>0$). This system
corresponds to a Hamiltonian system derived in \cite{CGK} (see their
formula (5.10)). As mentioned before, the present, three-parameter
family of systems allows one to circumvent the problem of
ill-posedness without the need of taking into account higher-order
terms in the expansion, as in \cite{CGK}).
\end{remark}
\begin{proof}  The proof is again made based on various
possibilities for the parameters in the problem. For this regime,
we have that $\eps\sim\mu\sim \eps_2\sim
\mu_2$ as $\eps\to 0$.  The overall idea of the argument is the same
as evinced in the proof of  Theorem \ref{propFDFD}.\\
 {\bf Step 1.}
{\emph The case $\alpha_1=0$, $\beta=0$, $\alpha_2=0$.}  Using
Remark \ref{remB} and Corollary \ref{coro2ter} (instead of
Proposition \ref{prop2} and Corollary \ref{coro2} as in the last
theorem) one checks immediately
 that
$$
    \nabla\psi_1=-\frac{1}{\gamma+\delta}\big[
    1+\mu\frac{1}{3\delta}\frac{1-\delta^2}{\gamma+\delta}\Delta
    +\eps_2\frac{1+\delta}{\gamma+\delta}\Pi(\zeta\cdot)\big]{\bf v}
    +O(\eps^2)
$$
(the nonlocal operator $\Pi$ does not appear in the final equations because
of the identity $\nabla\cdot \Pi V=\nabla\cdot V$ for all $V\in H^1(\R^d)^d$).\\
{\bf Step 2.} {\emph The case $\alpha_1 \geq 0$, $\beta=0$,
$\alpha_2=0$.} To use the BBM-trick, remark that for all
$\alpha_1\geq 0$,
$$
    \nabla\cdot{\bf v}=(1-\alpha_1)\nabla\cdot {\bf v}
    -\alpha_1(\gamma+\delta)\dt\zeta+O(\eps).
$$
Substitute  this relation into the third-derivative term of the
first
equation of the system derived in Step 1.\\
{\bf Step 3.} {\emph The case $\alpha_1\geq 0$, $\beta\geq 0$,
$\alpha_2=0$.} It suffices to replace ${\bf v}$ by
$(1-\mu\beta\Delta){\bf v}_\beta$ in the system of equations derived
in Step 2.\\
{\bf Step 4.} {\emph The case $\alpha_1\geq 0$, $\beta\geq 0$,
$\alpha_2\leq 1$.} This is exactly as in Step 4 of Theorem
\ref{propBFD}.
\end{proof}

\subsection{The Shallow Water/Shallow Water Regime:
$\mu\sim \mu_2\ll1$}\label{SW-SW} Contrary to the regimes
investigated above, large amplitude interfacial deformations are
allowed for both fluids, as $\eps\sim\eps_2=O(1)$. As in the
previous section, an asymptotic model can be derived from
(\ref{eqdepndim}) by replacing the operators $G^\mu[\eps\zeta]$ and
${\bf H}^{\mu,\delta}[\eps\zeta]$ by their asymptotic expansions,
provided by Proposition \ref{prop1} and Corollary \ref{coro1} in the
present regime. The following theorem shows that the internal wave
equations are consistent in this regime with the Shallow
water/Shallow water system,
\begin{equation}\label{eqSWSW}
    \left\lbrace
    \begin{array}{l}
    \dt \zeta+\frac{1}{\gamma+\delta}\nabla\cdot \big(h_1{\mathfrak Q}[\frac{\gamma-1}{\gamma+\delta}\eps\delta\zeta](h_2{\bf v})\big)=0,\vspace{1mm}\\
    \dt {\bf v}
    +(1-\gamma)\nabla\zeta\\
    \indent+\frac{\eps}{2}\nabla\Big(
	 \big\vert 
	{\bf v}-\frac{\gamma}{\gamma+\delta}{\mathfrak Q}[\frac{\gamma-1}{\gamma+\delta}\eps\delta\zeta](h_2{\bf v})\big\vert^2
    -\frac{\gamma}{(\gamma+\delta)^2} \big\vert
	{\mathfrak Q}[\frac{\gamma-1}{\gamma+\delta}\eps\delta\zeta](h_2{\bf v})\big\vert^2\Big)=0,
    \end{array}\right.
\end{equation}
where $h_1=1-\eps\zeta$, $h_2=1+\eps\delta\zeta$, and the operator ${\mathfrak Q}$ is defined in Lemma \ref{lemmproj}.
\begin{theorem}\label{propSWSW}
    Let $0<\delta^{min}<\delta^{max}\leq (1-\delta(1-H_1))^{-1}$.
    The internal waves equations
    (\ref{eqdepndim}) are consistent with the SW/SW equations
    (\ref{eqSWSW}) in the sense of Definition \ref{defcons},
    with a precision $O(\mu)$, and uniformly with respect
    to
    $\eps\in[0,1]$, $\mu\in (0,1)$ and
    $\delta\in [\delta^{min},\delta^{max}]$.
\end{theorem}
\begin{remark}
    Taking $\gamma=0$ and $\delta=1$
    in the SW/SW equations (\ref{eqSWSW}) yields
     the usual shallow water equations for surface water waves
	(recall that it follows from Lemma \ref{lemmproj} that
	$\nabla\cdot [(1-\eps\zeta){\mathfrak Q}[-\eps\zeta]
	((1+\eps\zeta){\bf v})]=\nabla\cdot ((1+\eps\zeta){\bf v})$).
\end{remark}
\begin{remark}
    In the one-dimensional case $d=1$, one has 
	$$
	\frac{1}{\gamma+\delta}{\mathfrak Q}
	[\frac{\gamma-1}{\gamma+\delta}\eps\delta\zeta](h_2{\bf v})
	=\frac{h_2}{\delta h_1+\gamma h_2}
	$$
    and the equations (\ref{eqSWSW}) take the simpler form
$$
  \left\lbrace
    \begin{array}{l}
    \dt \zeta+\partial_x\big(\frac{h_1h_2}{\delta h_1+\gamma  h_2}{\bf v}\big)
	=0,\vspace{1mm}\\
    \dt {\bf v}
    +(1-\gamma)\partial_x\zeta
    +\frac{\eps}{2}\partial_x\big(
	\frac{(\delta h_1)^2-\gamma h_2^2}{(\delta h_1+\gamma h_2)^2}
	\vert{\bf v}\vert^2\big)=0,
    \end{array}\right.
$$
	which coincides of course with the system
    (5.26) of \cite{CGK}. The presence of the nonlocal
    operator ${\mathfrak Q}$, which does not seem to have been noticed before,
    appears to be a purely two dimensional
    effect.
\end{remark}
\begin{proof}
First remark that with the range of parameters considered in the
theorem, one has $\mu\sim \mu_2$ as $\mu\to 0$ while
$\eps\sim \eps_2=O(1)$. \\
By the definition (\ref{defv}) of ${\bf v}$ and using Proposition \ref{prop1} and Corollary
\ref{coro1}, one deduces from (\ref{eqdepndim}) that
\begin{equation}\label{chicago}
    \left\lbrace
    \begin{array}{l}
    \dt \zeta-\nabla\cdot ((1-\eps\zeta)\nabla\psi_1)=O(\mu),\\
    \dt {\bf v}
    +(1-\gamma)\nabla\zeta
    +\frac{\eps}{2}\nabla(
    \vert {\bf H}^{\mu,\delta}[\eps\zeta]\psi_1\vert^2
    -\gamma\vert \nabla\psi_1\vert^2)=O(\mu).
    \end{array}\right.
\end{equation}
Recall now that ${\bf H}^{\mu,\delta}[\eps\zeta]\psi_1={\bf
v}+\gamma\nabla\psi_1$; since moreover one also gets from 
Corollary \ref{coro1}
that ${\bf H}^{\mu,\delta}[\eps\zeta]\psi_1=-\delta {\mathfrak Q}[\eps_2\zeta](h_1\nabla\psi_1)+O(\mu)$, it is straightforward to deduce that
$$
{\bf v}+\gamma\nabla\psi_1=-\delta {\mathfrak Q}[\eps_2\zeta](h_1\nabla\psi_1)+O(\mu)
$$
Multiplying this relation by $h_2$ and taking the divergence, one gets
\begin{eqnarray*}
	\nabla\cdot(h_2{\bf v})+\gamma\nabla\cdot(h_2\nabla\psi_1)&=&
	-\delta \nabla\cdot(h_2{\mathfrak Q}[\eps_2\zeta](h_1\nabla\psi_1))+\nabla\cdot O(\mu)\\
&=&-\delta \nabla\cdot(h_1\nabla\psi_1)+\nabla\cdot O(\mu),
\end{eqnarray*}
where the second equality comes from the definition of the operator ${\mathfrak Q}[\eps_2\zeta]$. We thus have
$$
\nabla\cdot((1+\frac{\gamma-1}{\gamma+\delta}\eps_2\zeta)\nabla\psi_1)=-\frac{1}{\gamma+\delta}\nabla\cdot(h_2{\bf v})+\nabla\cdot O(\mu),
$$
and we can therefore use Lemma \ref{lemmproj} to conclude that
$$
	\nabla\psi_1=-\frac{1}{\gamma+\delta}{\mathfrak Q}[\frac{\gamma-1}{\gamma+\delta}\eps_2\zeta](h_2{\bf v})+O(\mu)
$$
and consequently, 
\begin{eqnarray*}
	{\bf H}^{\mu,\delta}[\eps\zeta]\psi_1&=&{\bf v}+\gamma\nabla\psi_1\\
	&=&{\bf v}-\frac{\gamma}{\gamma+\delta}{\mathfrak Q}[\frac{\gamma-1}{\gamma+\delta}\eps_2\zeta](h_2{\bf v})+O(\mu).
\end{eqnarray*}
Replacing $\nabla\psi_1$ and ${\bf H}^{\mu,\delta}[\eps\zeta]\psi_1$ by these
two expressions in (\ref{chicago}) yields the result.
\end{proof}

\subsection{The Shallow Water/Small Amplitude Regime:
$\mu\ll1$, $\eps_2\ll1$}\label{SW-SA}

Derived here are various models corresponding to the case when the
upper fluid layer is shallow, but this restriction is not required
of the lower layer.  The interfacial deviations are thus not
necessarily small relative to the the upper fluid depth, but they
are small relative to the undistrubed depth of the lower layer.
Different systems of equations obtain, depending on the  sizes of
the paramaters $\eps$, $\mu$ and $\delta$ (and thus $\eps_2$ and
$\mu_2$).

\subsubsection{The Shallow Water/Full Dispersion Regime: $\mu\sim \eps_2^2\ll 1$, $\eps\sim\mu_2\sim 1$}\label{SW-FD}

In this regime, the internal waves equations are consistent with the
 Shallow Water/Full Dispersion system,
\begin{equation}\label{eqSWFD}
    \left\lbrace
    \begin{array}{l}
    \dsp \dt \zeta+\frac{1}{\gamma}\nabla\cdot (h_1{\bf v})
    -\frac{\sqrt{\mu}}{\gamma^2}
    \nabla\cdot\big(h_1\D\coth(\sqrt{\mu_2}\D)\Pi(h_1 {\bf
    v})\big)=0,
    \vspace{1mm}\\
    \dsp\dt {\bf v}
    +(1-\gamma)\nabla\zeta
    -\frac{\eps}{2\gamma}\nabla\big[
    \vert{\bf v}\vert^2-2\frac{\sqrt{\mu}}{\gamma}{\bf v}\cdot\big(\D\coth(\sqrt{\mu_2}\D)\Pi(h_1{\bf v})\big)\big]=0,
    \end{array}\right.
\end{equation}
where $h_1=1-\eps\zeta$ and $\Pi=-\frac{\nabla\nabla^T}{\Delta}$.
\begin{theorem}\label{propSWFD}
    Let $0<c^{min}<c^{max}$ and $\mu_2^{min}<\mu_2<\mu_2^{max}$.
    The internal waves equations
    (\ref{eqdepndim}) are consistent with the SW/FD equations
    (\ref{eqSWFD}) in the sense of Definition \ref{defcons},
    with a precision $O(\mu)$, and uniformly with respect
    to
    $\eps\in[0,1]$, $\mu\in (0,1)$ and
    $\delta\in (0,1)$ satisfying the conditions
    $$
    c^{min}<\frac{\mu}{\eps^2\delta^2}<c^{max}
    \quad\mbox{ and }\quad
    \mu_2^{min}<\frac{\mu}{\delta^2}<\mu_2^{max}.
    $$
\end{theorem}
\begin{remark}
    The SW/FD system (\ref{eqSWFD}), which as far as we know is new, is a generalization of
    the results of \S 5.4 of \cite{CGK} to the two-dimensional
    case $d=2$ and to the case of a lower layer of finite
    depth (the case of an infinite lower layer is formally recovered here
    by taking $\tmd=1$ in (\ref{eqSWFD})).
\end{remark}
\begin{proof}
First remark that with the range of parameters considered in the
theorem, one has $\eps_2^2\sim \mu$ and $\eps\sim\mu_2\sim 1$ as
$\mu\to 0$.\\ Proposition \ref{prop1} implies that
$\frac{1}{\mu}G^\mu[\eps\zeta]\psi_1=\nabla\cdot
(h_1\nabla\psi_1)+O(\mu)$ while it follows from the definition of
${\bf v}$ and Corollary \ref{coro3} that
$$
    \nabla \psi_1=-\frac{1}{\gamma}{\bf v}+\frac{\sqrt{\mu}}{\gamma^2}
    \frac{\D}{\tmd}\Pi(h_1{\bf v})+O(\mu).
$$
One then concludes the proof exactly as in the previous sections.
\end{proof}
\subsubsection{The Intermediate Long Wave Regime: $\mu\sim \eps^2\sim\eps_2\ll1$, $\mu_2\sim 1$}\label{ILW}

In this regime, a \emph{one-parameter family} of intermediate long
wave systems may be derived from the internal waves equations. These
depend upon the parameter $\alpha$ and have the form
\begin{equation}\label{eqILW}
    \left\lbrace
    \begin{array}{l}
    \dsp[1+\sqrt{\mu}\frac{\alpha}{\gamma}\D\coth(\sqrt{\mu_2}\D)]\dt \zeta+\frac{1}{\gamma}\nabla\cdot ((1-\eps\zeta){\bf v})\vspace{1mm}\\
    \indent\dsp-(1-\alpha)\frac{\sqrt{\mu}}{\gamma^2}
    \D\coth(\sqrt{\mu_2}\vert D\vert)\nabla\cdot{\bf v}=0\vspace{1mm},\\
    \dsp \dt {\bf v}
    +(1-\gamma)\nabla\zeta
    -\frac{\eps}{2\gamma}\nabla
    \vert{\bf v}\vert^2=0.
    \end{array}\right.
\end{equation}
\begin{theorem}\label{propILW}
    Let $0<c^{min}<c^{max}$, $\mu_2^{min}<\mu_2<\mu_2^{max}$.
    The internal wave equations
    (\ref{eqdepndim}) are consistent with the ILW system
    (\ref{eqILW}) in the sense of Definition \ref{defcons},
    with a precision $O(\mu)$, and uniformly with respect
    to
    $\eps\in[0,1]$, $\mu\in (0,1)$ and
    $\delta\in (0,1)$ satisfying the conditions
    $$
    c^{min}<\frac{\mu}{\eps^2}<c^{max}
    \quad\mbox{ and }\quad
    \mu_2^{min}<\frac{\mu}{\delta^2}<\mu_2^{max}.
    $$
\end{theorem}
\begin{remark}
In dimension $d=1$ and with $\alpha=0$, (\ref{eqILW}) corresponds to
(5.47) of \cite{CGK}. However this system is not linearly
well-posed.  It is straightforward to ascertain that the condition
$\alpha\geq 1$ insures that (\ref{eqILW}) is linearly well-posed for
either $d = 1$ or $d=2$.
\end{remark}
\begin{remark}
    The ILW equation derived in \cite{Joseph,KK} is obtained as the
    unidirectional limit of the one dimensional
    ($d=1$) version of (\ref{eqILW}) -- see for instance
    \S 5.5 of \cite{CGK}.
\end{remark}
\begin{proof}
{\bf Step 1.}  {\emph The case $\alpha = 0$.}   We are working with
the regime $\mu\sim \eps^2\sim\eps_2\ll1$ and $\mu_2\sim 1$ as
$\mu\to 0$. In this situation, Proposition \ref{prop1} allows us to
write $\frac{1}{\mu}G^\mu[\eps\zeta]\psi_1=\nabla\cdot
((1-\eps\zeta)\nabla\psi_1)+O(\mu)$ while it follows from the
definition of ${\bf v}$ and (\ref{asympILW}) that
$$
    \nabla \psi_1=-\frac{1}{\gamma}{\bf v}+\frac{\sqrt{\mu}}{\gamma^2}
    \frac{\D}{\tmd}{\bf v}+O(\mu).
$$
 Substituting these two relations into the internal wave equations
(\ref{eqdepndim}) leads
to the advertised result with $\alpha = 0$.\\
{\bf Step 2.} {\emph The case $\alpha\geq 0$.} This result follows
from Step 1 and the observation that
$$
    \nabla\cdot{\bf v}=(1-\alpha)\nabla\cdot {\bf v}-\alpha \gamma\dt \zeta+O(\eps,\sqrt{\mu}).
$$
As mentioned already, the restriction on $\alpha$ is not to obtain
consistency, but rather to ensure linear well-posedness.
\end{proof}

\subsubsection{The Benjamin-Ono Regime: $\mu\sim \eps^2\ll1$, $\mu_2=\infty$}\label{BO}

For completeness, we investigate the Benjamin-Ono regime,
characterized by the asumption $\delta=0$ (the lower layer is of
infinite depth). Taking $\mu_2=\infty$ in (\ref{eqILW}) leads one to
replace $\coth(\sqrt{\mu_2}\D)$ by $1$.   The following
two-dimensional generalization of the system (5.31) in \cite{CGK}
emerges in this situation.
\begin{equation}\label{eqBOsyst}
    \left\lbrace
    \begin{array}{l}
    [1+\sqrt{\mu}\frac{\alpha}{\gamma}\D]\dt \zeta+\frac{1}{\gamma}\nabla\cdot ((1-\eps\zeta){\bf v})
    -(1-\alpha)\frac{\sqrt{\mu}}{\gamma^2}
    \D\nabla\cdot{\bf v}=0,\\
    \dt {\bf v}
    +(1-\gamma)\nabla\zeta
    -\frac{\eps}{2\gamma}\nabla
    \vert{\bf v}\vert^2=0.
    \end{array}\right.
\end{equation}
Neglecting the $O(\sqrt{\mu})=O(\eps)$ terms, one finds that $\zeta$
must solve a wave equation (with speed
$\sqrt{\frac{1-\gamma}{\gamma}}$). Thus, in the case of horizontal
dimension $d=1$, any interfacial perturbation splits up at first
approximation into two counter-propagating waves. If one includes
the $O(\sqrt{\mu},\eps)$ terms, one obtains the one-parameter family
\begin{equation}\label{eqBO}
    (1+\sqrt{\mu}\frac{\alpha}{\gamma}\vert\partial_x\vert)\dt \zeta
    +c\partial_x\zeta-\eps \frac{3}{4}c\partial_x \zeta^2
    -\frac{\sqrt{\mu}}{2\gamma}c(1-2\alpha)
    \vert\partial_x\vert\partial_x\zeta=0,
\end{equation}
of \emph{regularized Benjamin-Ono equations} (see \ref{BK}).  Here,
$c=\sqrt{\frac{1-\gamma}{\gamma}}$. The usual Benjamin-Ono equation
is recovered by taking $\alpha=0$.

\appendix
\section{Proof of Proposition \ref{expgal}}\label{appA}

The proof is made in five steps.\\
{\bf Step 1.} \emph{Coercivity of the operator
$\nabla^{\mu_2}_{X,z}\cdot
Q^{\mu_2}[\eps_2\zeta]\nabla^{\mu_2}_{X,z}$.} Exactly as in Prop.
2.3 of \cite{AL}, one may check that
$$
    \forall \Theta\in \R^{d+1},\qquad
    \Theta\cdot Q^{\mu_2}[\eps_2\zeta]\Theta\geq \frac{1}{k}\vert\Theta\vert^2,
$$
with $k=k(\frac{1}{H_2},{\eps\sqrt{\mu}},\eps_2\vert \zeta
\vert_{W^{1,\infty}})>0$.\\
{\bf Step 2.} \emph{Existence of a unique solution  to
(\ref{eqlm})}. Owing to Step 1, existence of a solution and
uniqueness up to a constant is provided by classical theorems (e.g.
Section V.7 of \cite{Taylor1}), provided that the source terms and
Neumann conditions satisfy the compatibility condition
$$
    \int_{{\mathcal S}}\nabla^{\mu_2}_{X,z}\cdot {\bf h}=
    \int_{\{z=0\}}(\sqrt{\mu_2}\nabla\cdot V+{\bf e_z}\cdot{\bf h})
    -
    \int_{\{z=-1\}}{\bf e_z}\cdot{\bf h}.
$$
This latter restriction is valid in the present circumstances on account of the divergence theorem.\\
{\bf Step 3.} \emph{$L^2$-estimate on $\nabla_{X,z}^{\mu_2}u$.}
Multiplying (\ref{eqlm}) by $u$, integrating by parts on both sides,
and using the Neumann conditions leads to
$$
    \int_{{\mathcal S}}\nabla_{X,z}^{\mu_2}u\cdot
    Q^{\mu_2}[\eps_2\zeta]\nabla_{X,z}^{\mu_2}u=
    -\int_{\{z=0\}}V\cdot \sqrt{\mu_2}\nabla u+
    \int_{\mathcal S}{\bf h}\cdot \nabla_{X,z}^{\mu_2}u.
$$
A direct consequence of Step 1 and the Cauchy-Schwarz inequality is
the inequality
$$
    \Vert \nabla^{\mu_2}_{X,z}u\Vert^2
    \leq k \big(\Vert {\bf h}\Vert\,\Vert \nm u\Vert
    +\vert V\vert_{H^{1/2}}\vert \sqrt{\mu_2}\nabla u\vert_{H^{-1/2}}).
$$
It follows from the trace theorem that
\begin{eqnarray*}
    \vert \sqrt{\mu_2}\nabla u\vert_{H^{-1/2}}&\leq& \cst\big(
    \Vert \sqrt{\mu_2}\nabla u\Vert+
    \Vert \Lambda^{-1}\sqrt{\mu_2}\dz \nabla u\Vert\big)\\
    &\leq& \cst\big(\Vert \nm u\Vert+\sqrt{\mu_2}\Vert \nm
    u\Vert\big).
\end{eqnarray*}
It is concluded that
$$
    \Vert \nm u\Vert\leq C(\frac{1}{H_2},\eps_2^{max},\mu_{2}^{max},\vert\zeta\vert_{W^{1,\infty}})\big(\Vert {\bf h}\Vert+\vert V\vert_{H^{1/2}}\big).
$$
{\bf Step 4.} \emph{$H^s$-estimate ($s\geq 0$) on
$\nabla_{X,z}^{\mu_2}u$.} Let $v=\Lambda^s u$.   Multiplying
(\ref{eqlm}) by $\Lambda^s$ on both sides, it results that $v$
solves the system
$$
    \left\lbrace
    \begin{array}{l}
    \nabla^{\mu_2}_{X,z}\cdot Q^{\mu_2}[\eps_2\zeta]\nabla^{\mu_2}_{X,z}
    v=\nabla_{X,z}^{\mu_2}\cdot \widetilde{{\bf h}}, \qquad\mbox{ in }{\mathcal S},
    \vspace{1mm}\\
    \partial_n v_{\vert_{z=0}}=\sqrt{\mu_2}\nabla\cdot \Lambda^s V+
    {\bf e_z}\cdot \widetilde{{\bf h}}_{\vert_{z=0}}
    ,\qquad
    \partial_n v_{\vert_{z=-1}}={\bf e_z}\cdot \widetilde{{\bf h}}_{\vert_{z=-1}},
    \end{array}\right.
$$
with $\widetilde{\bf h}=\Lambda^s{\bf
h}+[Q^{\mu_2}[\eps_2\zeta],\Lambda^s]\nabla_{X,z}^{\mu_2}u$. From
Step 3 and the definition of $v$,  it is thus deduced that $\Vert
\Lambda^s \nabla_{X,z}^{\mu_2}u\Vert$ is bounded from above by
$$
     C(\frac{1}{H_2},\eps_2^{max},
    \mu_{2}^{max},\vert\zeta\vert_{W^{1,\infty}})
    \big(\Vert\Lambda^s {\bf h}\Vert +\vert V\vert_{H^{s+1/2}}
    + \Vert [Q^{\mu_2}[\eps_2\zeta],\Lambda^s]\nabla_{X,z}^{\mu_2}u\Vert
    \big).
$$
Using the expression for $Q^{\mu_2}[\eps_2 \zeta]$ and the commutator
estimate
$$
\vert [\Lambda^s,f]g\vert_{2}\leq C \vert \nabla
f\vert_{H^{\max\{t_0,s-1\}}}\vert g\vert_{H^{s-1}},
$$
which holds for some constant $C$ which depends upon $s >
-\frac{d}2$ and $t_0>\frac{d}2$ (see Th. 6 of \cite{LJFA}), we
obtain
$$
     \Vert [Q^{\mu_2}[\eps_2\zeta],\Lambda^s]\nabla_{X,z}^{\mu_2}u\Vert
    \leq
    C(\frac{1}{H_2},\eps_2^{max},\mu_2^{max},\vert\zeta\vert_{H^{\max\{t_0+2,s+1\}}})
    \Vert \Lambda^{s-1}\nabla_{X,z}^{\mu_2}u\Vert
    \big).
$$
We thus get an estimate on $\Vert \Lambda^s
\nabla_{X,z}^{\mu_2}u\Vert$ in terms of $\Vert
\Lambda^{s-1}\nabla_{X,z}^{\mu_2}u\Vert$ which, together with Step 3
(i.e. $s=0$) allows us to derive the following relation by induction
(and interpolation when $s\in (0,1)$):
$$
    \forall s\geq 0,\qquad
    \Vert \Lambda^s \nabla_{X,z}^{\mu_2}u\Vert\leq
    C(\frac{1}{H_2},\eps_2^{max},
    \mu_2^{max},\vert\zeta\vert_{H^{\max\{t_0+2,s+1\}}})
    (\Vert \Lambda^s{\bf h}\Vert+\vert V\vert_{H^{s+1/2}}).
$$
{\bf Step 5.}  \emph{$H^s$-estimate ($s\geq 0$) on
$\dz\nabla_{X,z}^{\mu_2}u$.} First remark that using the equation
yields
\begin{eqnarray*}
    \lefteqn{\frac{1+\mu\eps^2(z+1)^2\vert\nabla\zeta\vert^2}{1+\eps_2\zeta}\dz^2u=
    \nm\cdot{\bf h}}\\
    & &-\sqrt{\mu_2}\nabla\cdot\big((1+\eps_2\zeta)\sqrt{\mu_2}\nabla u
    -\sqrt{\mu}\eps(z+1)\nabla\zeta\dz u\big)\\
    & &+\sqrt{\eps}\mu\nabla\zeta\cdot(\sqrt{\mu_2}\nabla u)
    -2\mu\eps^2(z+1)\frac{\vert\nabla\zeta\vert^2}{1+\eps_2\zeta}\dz u,
\end{eqnarray*}
from which one obtains the estimate
$$
    \Vert \Lambda^s\dz^2 u\Vert\leq
    C(\eps_2^{max},\mu_2^{max},\vert \zeta\vert_{H^{\max\{t_0+2,s+1\}}})
    (\Vert \Lambda^s\nm {\bf h}\Vert+\sqrt{\mu_2}\Vert \Lambda^{s+1}
    \nm u\Vert).
$$
Use this to write
\begin{eqnarray*}
    \Vert \Lambda^s\dz\nm u\Vert&\leq& \sqrt{\mu_2}\Vert\Lambda^s\dz \nabla u\Vert
    +\Vert \Lambda^s\dz^2 u\Vert\\
    &\leq&  C(\eps_2^{max},\mu_2^{max},\vert \zeta\vert_{H^{\max\{t_0+2,s+1\}}})
    (\Vert \Lambda^s\nm {\bf h}\Vert+\sqrt{\mu_2}\Vert \Lambda^{s+1}
    \nm u\Vert).
\end{eqnarray*}
With the help of Step 4, one obtains the inequality
$$
    \Vert \Lambda^s\dz\nm u\Vert
    \leq
    C(\frac{1}{H_2},\eps_2^{max},\mu_2^{max},\vert \zeta\vert_{H^{\max\{t_0+2,s+1\}}})
    (\Vert{\bf h}\Vert_{H^{s+1,1}}+\vert V\vert_{H^{s+3/2}}).
$$
{\bf Step 6.} \emph{Conclusion.} By the trace theorem we may assert
that for all $s\geq 0$,
$$
    \vert \nabla u_{\vert_{z=0}}\vert_{H^s}\leq
    \cst \Vert \nabla u\Vert_{H^{s+1/2,1}}
    \leq \frac{\cst}{\sqrt{\mu_2}}
    \Vert \nabla^{\mu_2}_{X,z} u\Vert_{H^{s+1/2,1}}.
$$
The desired result now follows from Steps 4 and 5.

\end{document}